\documentstyle[12pt]{article}
\newcommand{\bold}{\bf}
\newcommand{\rbox}[1]{{\rm{#1}}}
\newcommand{\bbox}[1]{{\bf{#1}}}          
\newcommand{\rmt}{\rm}

\newcommand{\bft}{\bf}
\newcommand{\bald}[1]{{#1}} 
\def\restriction{\mathbin{\hspace{0.1ex}|\hspace{0.1ex}}}
\newtheorem{theorem}{Theorem}
\newtheorem{assertion}[theorem]{Assertion}
\newtheorem{corollary}[theorem]{Corollary}
\newtheorem{definition}{Definition}

\newtheorem{lemma}[theorem]{Lemma}
\newtheorem{proposition}[theorem]{Proposition}
\newtheorem{remark}[theorem]{Remark}
\newtheorem{warning}[theorem]{Warning}
\newcommand{\proof}{\noi{\bft Proof\hspace{2mm} }}
\newcommand{\TF}{\it}
\newcommand{\bass}{\begin{assertion}\TF\ } 
\newcommand{\eass}{\end{assertion}} 
\newcommand{\bcor}{\begin{corollary}\TF\ }
\newcommand{\ecor}{\end{corollary}}
\newcommand{\bdf} {\begin{definition}\rmt}
\newcommand{\edf} {\end{definition}}
\newcommand{\ble} {\begin{lemma}\TF\ }
\newcommand{\ele} {\end{lemma}}
\newcommand{\bte} {\begin{theorem}\TF\ }
\newcommand{\ete} {\end{theorem}}
\newcommand{\bpro}{\begin{proposition}\TF\ } 
\newcommand{\epro}{\end{proposition}} 
\newcommand{\brem} {\begin{remark}\rmt\ }
\newcommand{\erem} {\end{remark}}
\newcommand{\bwa} {\begin{warning}\rmt\ }
\newcommand{\ewa} {\end{warning}}
\newcommand{\qed} {\hfill$\msur\Box\msur$} 
\newcommand{\ben}{\begin{enumerate}}
\newcommand{\een}{\end{enumerate}}
\newcommand{\bit}{\begin{itemize}}
\newcommand{\eit}{\end{itemize}}
\newcommand{\bce}{\begin{center}}
\newcommand{\ece}{\end{center}}
\newcommand{\bde}{\begin{description}}
\newcommand{\ede}{\end{description}}
\newcommand{\bay}{\begin{array}}
\newcommand{\eay}{\end{array}}
\newcommand{\bqu}{\begin{quotation}}
\newcommand{\equ}{\end{quotation}}

\newcommand{\ZFC} {\bbox{ZFC}}
\newcommand{\IS}  {\bbox{IS}}
\newcommand{\cs}  {\bbox{CS}}
\newcommand{\bbsp}{\hspace{0.5pt}}
\def\a{{\bbsp\bbox{a}\bbsp}}
\def\x{{\bbsp\bbox{x}\bbsp}}
\def\y{{\bbsp\bbox{y}\bbsp}}
\def\i{\bald{i}}
\def\j{\bald{j}}
\newcommand{\dom}{\rbox{dom}\,}

\newcommand{\card}{\rbox{card}\,}
\newcommand{\ste}{\rbox{stem}\hspace{1pt}}
\newcommand{\roo}{\rbox{root}\hspace{1pt}}
\newcommand{\Ord}{\rbox{Ord}}

\newcommand{\perfm}{{\dP}}
    \newcommand{\peM}[1]{\perfm_{#1}}
       \newcommand{\peim}{\perfm}
\newcommand{\contm}{{\dF}}  
    \newcommand{\cntm}[1]{\contm_{#1}}
        \newcommand{\cntim}{\contm}
\newcommand{\perf}{\rbox{Perf\hspace{0.5pt}}}
    \newcommand{\pe}[1]{\perf_{#1}}
        \newcommand{\pei}{\perf}
    \newcommand{\ipe}[1]{\perf'_{#1}}
        \newcommand{\ipei}{\perf'}
\newcommand{\cont}{\rbox{Cont\hspace{0.5pt}}}
    \newcommand{\cnt}[1]{\cont_{#1}}
        \newcommand{\cnti}{\cont}
\newcommand{\can}[1]{\cD^{#1}}
        \newcommand{\cani}{\can{\bI}}
\newcommand{\diam}{\rbox{diam}\hspace{1pt}}
\newcommand{\spl}{\rbox{Spl}}
\newcommand{\clo}[1]{\rbox{Clo}\hspace{0.5pt}[#1]}
\newcommand{\rL}{\rbox{L}}
\newcommand{\rV}{\rbox{V}}

\newcommand{\ima}{{\hbox{\hspace{1pt}\rmt ''}}}
\newcommand{\al}{\alpha} 
\newcommand{\da}{\delta}
\newcommand{\ga}{\gamma}
\newcommand{\La}{\Lambda} 
\newcommand{\la}{\lambda} 
\newcommand{\vpi}{\varphi}
\newcommand{\om}{\omega} 
\newcommand{\Om}{\Omega} 
\newcommand{\vep}{\varepsilon}
\newcommand{\za}{\zeta}
\newcommand{\fs}[2]{{\bold\Sigma}^{#1}_{#2}}

\newcommand{\iSigma}{{\mathchar"7106}}
\newcommand{\is}[2]{\iSigma^{#1}_{#2}}

\newcommand{\bbb}{\hspace{0.5pt}} 

\newcommand{\dP}{\mathord{{\rm I}\hspace{-2.5pt}{\rm P}}}
\newcommand{\dZ}{\mathord{{\sf Z}\hspace{-4.5pt}{\sf Z}}}
\newcommand{\dF}{\mathord{{\rm I}\hspace{-2.5pt}{\rm F}}}
\newcommand{\bI} {\bbox{I}} 
\newcommand{\bJ} {\bald{J}} 
\newcommand{\gM}{{\bbb{\bf M}\bbb}}
\newcommand{\gN}{{\bbb{\bf N}\bbb}}
\newcommand{\skrsp}{\hspace{0.5pt}}
\newcommand{\skri}[1]{{\skrsp {\cal #1}\skrsp}}
\newcommand{\cD}{{\skri D}}
\newcommand{\cN}{{\skri N}}

\newcommand{\cX}{{\skri X}}
\newcommand{\sq}{\subseteq}
\newcommand{\cj}{\;\,\&\;\,}

\newcommand{\lra}{\longrightarrow} 
\newcommand{\llra}{\longleftrightarrow} 
\newcommand{\res}{{\hspace{1.5pt}\restriction\hspace{1.5pt}}}
    \newcommand{\rsd}[1]{\mathord{\restriction_{#1}}}
\newcommand{\ares}{\hspace{-2pt}\restriction^{-1}\hspace{-1pt}}
\newcommand{\we}{{\mathbin{\kern 1.3pt ^\wedge}}}
\newcommand{\<}{\leq}
\newcommand{\lsp}{\hspace{-2pt}}
     \def\>{\geq}
\newcommand{\ti}{\times}
\newcommand{\dm}{$$}
\newcommand{\ang} [1]{\langle #1\rangle}
\newcommand{\ans} [1]{\{\hspace{0.2mm}#1\hspace{0.2mm}\}}
\newcommand{\dd}[1]{$\kern-0.7mm{#1}\kern-1mm$-}
\newcommand{\noi}{\noindent}

\newcommand{\vom}{\vspace{1mm}}
\newcommand{\its}{\vspace{-1mm}}
\newcommand{\unf}{{\underline f}}
\newcommand{\unc}{{\underline c}}

\newcommand{\supp}[1]{\|#1\|}
\font\ess=cmss9
\scriptfont\sffam=\ess
\newcommand{\refo}[1]{{\mathchoice {\hbox{\sf{#1}}} 
{\hbox{\sf{#1}}}{\hbox{\ess{#1}}}{\hbox{\ess{#1}}} }}
\newcommand{\relf}[1]{{\refo #1}}
\def\C  {\mathbin{\relf{C}}}

\def\E  {\mathbin{\relf{E}}}
\def\Eo {\mathbin{\relf{E}_0}}

\newcommand{\msur}{\hspace{-1\mathsurround}}

\newcommand{\itla}{\item\label}
\begin{document}

\title{On non--wellfounded iterations of perfect set forcing with 
application to the Glimm -- Effros property} 

\author{Vladimir Kanovei
\thanks{Moscow Transport Engineering Institute}
\thanks{{\tt kanovei@nw.math.msu.su} \ and \ 
{\tt kanovei@math.uni-wuppertal.de}
}
\thanks{Partially supported by AMS grant} 
}
\date{01 August 1995} 
\maketitle
\normalsize


\begin{abstract}\vspace{0mm} 
\noi
We prove that if $\bI$ is a p.\ o. set in a countable 
transitive model $\gM$ of $\ZFC$ then $\gM$ can be extended by a 
generic sequence of reals $\a_\i,$ $\i\in\bI,$ such that 
$\aleph_1^\gM$ is preserved and every $\a_\i$ is Sacks generic 
over $\gM[\ang{\a_\j:\j<\i}]$.

The structure of the degrees of \dd\gM constructibility of reals 
in the extension is investigated.

As an application, we obtain a model in which the $\is12$ 
equivalence relation $x \E y$ iff $\rL[x]=\rL[y]$ ($x,\,y$ are 
reals) does not admit a reasonable form of the Glimm -- Effros 
theorem.

\vfill

\vfill

\bce
{\bft Acknowledgements}
\ece
\noi
The author is in debt to M.~J.~A.~Larijani, the president of 
IPM (Tehran, Iran), for the support during the initial period of 
work on this paper in May and June 1995.
The author is pleased to thank S.\ D.\ Friedman, M.\ Groszek, 
G.\ Hjorth, \hbox{A.\ S.\ Kechris,} A.\ W.\ Miller, and T.\ Slaman 
for useful discussions and interesting information on Solovay model, 
and the Glimm -- Effros matter, and iterated Sacks forcing.
\vfill

\end{abstract}

\vfill

\subsection*{Introduction}

It is a common practice in set theory that one is interested to 
consider a generic extension $M_1$ of a model $M,$ after this a 
generic extension $M_2$ of $M_1,$ and so on, including the case 
of infinite or transfinite number of steps. Iterated forcing of 
Solovay and Tennenbaum~\cite{st} allows us to engineer this 
iterated construction in an ordinary one--step generic extension. 

In the most of cases iterated forcing is used to define 
transfinite sequences of models such that every model is a certain 
generic extension of the preceding model. (We do not consider here  
sophisticated details at limit steps). Identifying the steps of 
this construction with ordinals, and interpreting the set of the 
ordinals involved as the {\it support\/} or the {\it ``length''\/} of 
the iteration, we may say that the classical iterated forcing is 
an iterated forcing of {\it wellordered\/} ``length''. 

In principle it does not require an essential improvement of the 
basic iterater forcing method to define iterations of 
{\it wellfounded\/}, but not linearly ordered, ``length''. 
This version is much rarely used then 
the basic one. (See Groszek and Jech~\cite{gj} for several 
known applications.) 

It is a much more challenging question (we refer to Groszek and 
Jech~\cite{gj}, p.~6) to carry out ``ill''founded iterations. No 
general method is known, at least. 

For a few number of rather simple forcing notions, ``ill''founded 
iterations can be obtained without any use of the idea of 
iteration at all. For example if $a\in 2^\om$ is a Cohen generic 
real over a model $M$ and $b=o(a)\in 2^m$ is defined for any 
$a\in 2^\om$ by $b(m)=a(2m),\;\forall\,m,$ then the sequence 
of reals $a_n$ defined by $a_0=a$ and $a_{n+1}=o(a_n)$ realizes 
the iteration of Cohen forcing of ``length'' $\om^\ast$ (the reverse 
order on natural numbers): every $a_n$ is Cohen generic over 
$M[\ang{a_m:m>n}].$ This construction can be applied to Solovay 
random reals as well.

An idea how to carry out iterated forcing of a linear but not 
wellordered ``length'' $I$ can be as follows. let us first consider a 
usual iteration of a ``length'' $\la\in\Ord$ as a pattern to follow. 
The forcing conditions in this case are functions $p$ defined on 
the set $\la=\ans{\al:\al<\la}$ and satisfying certain property 
$P(p,\al)$ 
for every $\al<\la.$ Now to proceed with the \dd\bI case one may 
want to use functions $p$ defined on $I$ and satisfying $P(p,i)$ 
for all $i\in I$. 

The principal problem in this argument is that in the wellordered 
setting the property $P(p,\al)$ is itself defined by induction on 
$\al$ in a quite sophisticated way. So we first have to eliminate 
the induction and extend the property $P$ to ``ill''ordered sets. 

We do not know how this can be realized at least for a more or 
less representative category of forcing notions. There is, 
however, a forcing which allows to express the property 
$P$ in simple geometrical terms, so that the ``ill''founded 
iterations become available. This is the perfect set forcing 
introduced by 
Sacks~\cite{sa}. (We refer to Baumgartner and Laver~\cite{bl} on 
matters of iterated Sacks forcing, and Groszek~\cite{g} on 
further applications.) 
\pagebreak[3]

\bte
\label{m}
Let\/ $\gM$ be a countable transitive model of $\ZFC,$ $\bI$ a 
partially ordered\nopagebreak{} 
set in\/ $\gM.$ Then there exists a generic\/ 
\dd{\aleph_1}preserving extension\/ 
$\gN=\gM[\ang{\a_\i:\i\in\bI}]$ of\/ $\gM$ such that\its
\ben
\def\theenumi{{\rmt\arabic{enumi}}}
\itla{m1} 
For every\/ $\i\in \bI,$ $\a_\i$ is a Sacks--generic real
over\/ $\gM[\ang{\a_\j:\j<\i}]$.\its

\itla{m1+} 
If\/ $\i,\,\j\in\bI$ and\/ $\i<\j$ then\/ $\a_\i\in\gM[\a_\j]$.\its

\itla{m2} 
If\/ $\xi\in\gM$ is an initial segment in\/ $\bI$ and\/ 
$\i\in\bI\setminus\xi$ then\/ 
$\a_\i\not\in\gM[\ang{\a_\j:\j\in\xi}]$.\its

\itla{m2+} 
If\/ $\xi\in\gM$ is a countable in\/ $\gM$ initial segment in\/ 
$\bI$ and\/ 
$c$ is a real in\/ $\gN$ such that\/ $\a_\i\in\gM[c]$ for all\/ 
$\i\in\xi$ then the indexed set\/ $\ang{\a_\i:\i\in\xi}$ 
belongs to\/ $\gM[c]$.\its

\itla{m3} 
For any initial segment\/ $\xi\sq\bI,$ $\xi\in\gM,$ and any real\/ 
$c\in\gN,$ we have exactly one from the following:\hspace{4mm} 
\begin{minipage}[t]{0.6\textwidth}
{\rmt(a)} $c\in \gM[\ang{\a_\i:\i\in\xi}]$, \hspace{10mm} 
or\vspace{2mm}

{\rmt(b)} \hspace{1mm}there exists\/ $\i\in\bI\setminus\xi$ 
such that $\a_\i\in \gM[c]$.
\end{minipage}
\hfill $\,$
\een
\ete
The set $\bI$ is not necessarily wellfounded or linearly ordered 
in $\gM.$ In the particular case of inverse ordinals taken as 
$\bI$ the theorem was recently proved by Groszek~\cite{g94}.

Items \ref{m2}, \ref{m2+}, \ref{m3} seem to show that the 
\dd\gM degree of a real $c$ in the extension $\gN$ intends to be 
determined by the set $\bI_c=\ans{\i\in\bI:\a_\i\in\gM[c]},$ an 
initial segment of $\bI$ by item~\ref{m1+}. One can easily prove 
that in fact $\bI_c=\bI_{c'}$ implies $\gM[c]=\gM[c']$ provided 
the set $\bI_c=\bI_{c'}$ belong to $\gM$ (e.g. in the case when 
all initial segments of $\bI$ belongs to $\gM$), but the general 
case remains open. 

The proof of the theorem is based on a version of iterated Sacks 
forcing realized in the form of perfect sets~\footnote
{\rmt\ We consider perfect sets rather than perfect trees, because 
the particular combinatorial properties we need hardly can be 
expressed in a reasonable form for the forcing realized using 
trees rather than perfect sets. Of course the absoluteness of the 
conditions is lost because a perfect set in $\gM$ is not perfect 
both in the universe and the extension, but this is a comparably 
minor problem, easily fixed by taking the topological closure.} 
with certain combinatorial properties.

We shall be mostly concentrated on the case when $\bI$ is 
{\it finite or countable\/} in $\gM.$ The forcing we use in this 
case will be a collection $\peim$ of perfect subsets of the 
product in $\gM$ of $\bI$ copies of the Cantor space $2^\om,$ 
\dd\i th copy being 
responsible for the corresponding $\a_\i.$ Sections \ref{prelim} 
through \ref{cont} of the paper present useful properties of sets 
in $\peim$ and continuous functions defined on them. This part of 
the paper is not related to any particular model but 
finally the reasoning will be ``relativized'' to $\gM$.

The results of this study are used in sections \ref{re} through 
\ref{uncount} for the proof of Theorem~\ref{m}. We show that the 
reals in a \dd\dP generic extension can be presented by continuous 
functions in the ground model $\gM,$ defined on sets $X\in\peim.$ 
It occurs that notions related to degrees of \dd\gM 
constructibility of reals in the extension are adequately 
reflected in properties of continuous functions in the ground 
model. 

The case of {\it arbitrary\/} $\bI$ is reduced in 
Section~\ref{uncount} to the case of countable $\bI$ by an 
ordinary ``countable support'' argument.

\subsubsection*{An application: non--Glimm--Effros $\is12$ 
equivalence}

Harrington, Kechris, and Louveau~\cite{hkl} proved that each 
Borel equivalence relation $\E$ on reals satisfies one and 
only one of the following conditions:\its
\ben
\def\theenumi{{\rmt\hskip2pt(\Roman{enumi})\hskip2pt}}
\def\labelenumi{\theenumi}
\itla{1} \msur 
$\E$ admits a countable Borel separating family.\its

\itla{2} \msur 
$\E$ continuously embeds $\Eo,$ the Vitali equivalence.
\een
(Some notation. A {\it separating family\/} for an equivalence 
$\E$ on reals is an indexed family $\ang{X_\al:\al<\ga}$ 
($\ga\in\Ord$)  
of sets $X_\al$ such that $x\E y$ iff 
$\forall\,\al\,(x\in X_\al\,\llra\,y\in X_\al)$ for all $x,\,y.$ 
$\Eo$ is the {\it Vitali equivalence\/} on the 
{\it Cantor space\/} $\cD=2^\om,$ defined 
by: ${x\Eo y}$ iff $x(n)=y(n)$ for all but finite $n\in\om.$ An 
{\it embedding\/} of $\Eo$ into $\E$ is a $1-1$ function 
$U:\cD\,\lra\,\hbox{reals}$ such that 
$x\Eo y\;\llra\;U(x)\E U(y)$ for all $x,\,y\in \cD$. 
We refer the reader to \cite{hkl} as the basic sourse of 
information on the matter.) 

Hjorth and Kechris~\cite{hk}, Hjorth~\cite{h-det}, 
Kanovei~\cite{k-sm,k-s11} obtained partial results of 
this type for $\fs11$ and even more complicated relations, which 
we do not intend to discuss here. 

However there exists a $\is12$ equivalence relation which does 
not admit a theorem of the Glimm--Effros type in $\ZFC,$ at least 
in the field of real--ordinal definable (R-OD, in brief) 
separating families and embeddings. 

\bte
\label{ge}
It is consistent with\/ $\ZFC$ that the\/ $\is12$ equivalence 
relation\/ $\C$ defined on reals by\/ $x\C y$ iff\/ 
$\rL[x]=\rL[y]\;:$\its
\bit
\item[--] neither has a R-OD separating family$;$
\its

\item[--] nor admits an uncountable R-OD pairwise\/ 
\dd\C inequivalent set$.$
\eit
\ete

\noi
{\bf Remarks} \ 

1. The ``nor'' part of the theorem implies that 
$\C$ does not embed $\Eo$ via a \hbox{R-OD} embedding, because 
obviously there exists a perfect set of pairwise \dd\Eo 
inequivalent points. 

2. It makes no sense to look for non-R-OD separating families 
in the ``either'' part. Indeed let $\kappa$ be 
the cardinal of the quotient set $\hbox{reals}/\C.$ Then any 
enumeration $\ang{X_\al:\al<\kappa}$ of all \dd\E equivalence 
classes is a separating family, but this construction does not 
guarantee the real--ordinal definability of the enumeration even 
in the case when $\E$ itself is R-OD (take $\Eo$ as an example). 

\vspace{4mm}  

The model for Theorem~\ref{ge} we propose is the iterated Sacks 
extension of the constructible model having $\om_1\ti\dZ$ 
($\om_1$ copies of the integers) as the ``length'' of iteration. 
The model is considered in Section~\ref{seqge}. 

\newpage
\subsection{Notation and pre--conditions}
\label{prelim}
\label{prop}
\label{splt}

{\it The ``length''\/}. Let $\bI$ be a fixed countable partially 
ordered set, which will be the ``length'' of iteration. Characters 
$\i,\,\j$ are used to denote elements of $\bI.$ Subsets of $\bI$ 
will be denoted by Greek letters $\xi,\,\eta,\,\za$.\vom\vom

{\it Spaces\/}. $\cN=\om^\om$ is the {\it Baire space\/}; points of 
$\cN$ will be called {\it reals\/}. $\cD=2^\om$ is the 
{\it Cantor space\/}. For 
$\xi\sq \bI,$ $\can\xi$ is the product of \dd\bI many copies of 
$\cD$ with the product topology (here $\bI$ is considered as 
discrete). Then every $\can\xi$ is a compact space homeomorphic 
to $\cD$ itself unless $\xi=\emptyset$.\vom\vom

{\it Projections\/}. Assume that $\xi\sq\eta\sq\bI.$ If 
$x\in\can\eta$ then let $x\res\xi\in\can\xi$ denote the 
usual restriction. If $X\sq\can\eta$ then let 
$X\res \xi=\ans{x\res \xi:x\in X}$.

But if $X\sq\can\xi$ then we set 
$X\ares \eta=\ans{y\in\can\eta:y\res \xi\in X}$.

In addition, if $\i\in\xi\sq\bI$ and $X\sq\can\xi$ then we put
$X(\i)=\ans{x(\i):x\in X}$.\vom\vom

{\it Initial segments\/}. Let $\IS$ denote the set of all 
initial segments of $\bI.$ 
For any $\i\in\bI,$ we put $[<\lsp\i]=\ans{\j\in\bI:\j<\i},$  
$[\not>\lsp\i]=\ans{\j\in\bI:\j\not>\i},$ and $[\<\lsp\i],$ 
$[\not\>\lsp\i]$ in the same way. 
To save space, let $X\rsd{<\i}$ mean $X\res[<\lsp\i],$ 
$\can{\<\i}$ mean $\can{[\<\i]},$ etc.\vom\vom

If $\i\in\xi\in\IS$ and $X\sq\can\xi,$ then we 
define $D_{Xz}(\i)=\ans{x(\i):x\in X\cj z=x\rsd{<\i}}$ 
for every $z\in X\rsd{<\i}.$ Thus $D_{Xz}(\i)\sq \cD$.\vom\vom

{\it Pre--conditions\/}.
The following definition would be sufficient for the purpose to 
prove Theorem~\ref{m} at least in two particular cases: when $\bI$ 
is wellfounded, and when $\bI$ is linearly ordered. In fact we 
don't know whether it gives the expected result in general case. 
We are not able to prove a very important technical fact 
(Proposition~\ref{clop} below): if $X'$ is a clopen (in the 
relative topology) nonempty subset of $X\in\ipe\xi$ then $X'$ 
contains a subset $X''\in\ipe\xi.$ This is why one more 
requirement will be added in Section~\ref{forc}, to define the 
notion of forcing completely. 

\bdf
(Pre--conditions)\\[1pt]
For any $\za\in\IS,$ $\ipe\za$ is the collection of 
all sets $X\sq\can\za$ such that\its
\ben
\def\theenumi{{\rmt P-\arabic{enumi}}}
\item The set $X$ is closed and nonempty.\its

\itla{perf1} 
If $\i\in\za$ and $z\in X\rsd{<\i}$ then 
$D_{Xz}(\i)$ is a perfect set in $\cD$.\its

\itla{oz} 
If $\i\in \za$ and $G\sq\cD$ is open then the set 
$\ans{x\rsd{<\i}:x\in X\cj x(\i)\in G}$ is open in 
$X\rsd{<\i}$.~\footnote
{\rmt\ In other words, it is required that the projection from 
$X\rsd{\<\i}$ to $X\rsd{<\i}$ is an open map.}\its

\itla{indep} 
If $\xi,\,\eta\in\IS,$ 
$\xi\cup \eta\sq \za,$ $x\in X\res \xi,$ $y\in X\res \eta,$ 
and $x\res (\xi\cap \eta)=y\res (\xi\cap \eta),$ then 
$x\cup y\in X\res (\xi\cup \eta)$.\its
\een
Finally we set $\ipei=\ipe\bI$.\qed
\edf

This section contains several quite elementary lemmas on 
pre--conditions.

\bass
\label{less'}
If\/ $X\in\ipe\za$ and\/ $\xi\in\IS,$ $\xi\sq \za,$ then\/ 
$X\res \xi\in\ipe\xi$.\qed
\eass

\ble 
\label{pro'}
Suppose that\/ $\xi,\,\eta,\,\za\in\IS,$ $\xi\cup \eta\sq \za,$ 
$X\in\ipe\za,$ $Y\sq X\res \eta,$ and\/ $Z=X\cap (Y\ares\za).$ 
Then $Z\res \xi=(X\res \xi)\cap (Y\res(\xi\cap \eta)\ares \xi)$.
\ele
\proof The inclusion $\sq$ is quite easy. To prove the opposite 
direction let $x$ belong to the right--hand side. Then in 
particular $x\res(\xi\cap \eta)=y\res(\xi\cap \eta)$ for some 
$y\in Y.$ On the other hand $x\in X\res\xi$ and $y\in X\res\eta.$ 
Condition~\ref{indep} implies $x\cup y\in X\res(\xi\cup \eta).$ 
Therefore $x\cup y\in Z\res(\xi\cup \eta)$ because $y\in Y.$ We 
conclude that $x\in Z\res \xi$.\qed

\ble
\label{apro'}
Suppose that\/ $\xi,\,\za\in\IS,$ $\xi\sq \za,$ $X\in\ipe\za,$ 
$Y\in\ipe\xi,$ and\/ $Y\sq X\res \xi.$ Then\/ 
$Z=X\cap (Y\ares \za)$ belongs to $\ipe\za$. 
\ele
\proof We check condition~\ref{perf1}. Let $\i\in\za$ and 
$z\in Z\rsd{<\i}.$ If $\i\in \xi$ then obviously 
$D_{Zz}(\i)=D_{Yz}(\i).$ If $\i\in \za\setminus\xi$ then 
$D_{Zz}(\i)=D_{Xz}(\i)$ by Lemma~\ref{pro'} (for 
$\eta=[\<\lsp\i]$).\vom

We check \ref{oz}. Let $\i\in\za.$ The case $\i\in\xi$ is easy as 
above, so let us suppose that $\i\in \za\setminus \xi.$ We assert 
that\its
\ben
\def\theenumi{(\arabic{enumi})}
\def\labelenumi{\theenumi}
\itla{aa1} \mathsurround=1mm
$Z\rsd{<\i}=(X\rsd{<\i})\cap (Y'\ares[<\lsp\i]),$ where 
$Y'=Y\res \xi'$ and $\xi'=\xi\cap[<\lsp\i]$, \hfill and\its

\itla{aa2} \mathsurround=1mm
$\ans{z\in Z:z(\i)\in G}\rsd{<\i}=
(\ans{x\in X:x(\i)\in G}\rsd{<\i})\cap (Y'\ares [<\lsp\i])$.\its
\een
Indeed \ref{aa1} immediately follows from Lemma~\ref{pro'}. 
The direction $\sq$ in \ref{aa2} is obvious. To prove the 
opposite direction, let $z\in\can{<\i}$ belong to the right--hand 
side, so that $z=x\rsd{<\i}$ for some $x\in X$ such that 
$x(\i)\in G,$ and $z\res\xi'=y\res\xi'$ for some $y\in Y.$ 
Applying property \ref{indep} of $X,$ we get $x'\in X$ such that 
$x'\rsd{\<\i}=x\rsd{\<\i}$ and $x'\res\xi=y.$ In particular, we 
have $x'\in Z$ and $x'(\i)\in G,$ so that 
$z=x\rsd{<\i}=x'\rsd{<\i}$ belongs to the left--hand side.

Thus both \ref{aa1} and \ref{aa2} are verified. Now it 
suffices to recall that $X$ satisfies~\ref{oz}, which implies that 
the set $\ans{x\in X:x(\i)\in G}\rsd{<\i}$ is clopen in 
$X\rsd{<\i}$.\vom

We check \ref{indep}. Assume that $\eta,\,\tau\in\IS,$ 
$x\res\eta\in Z\res\eta,$ and $x\res\tau\in Z\res\tau;$ we have 
to prove that $x\res(\eta\cup\tau)\in Z\res(\eta\cup\tau).$ Let 
$\eta'=\xi\cap\eta$ and $\tau'=\xi\cup\tau$.

First of all we note that 
$x\res(\eta\cup\tau)\in X\res(\eta\cup\tau)$ by 
property~\ref{indep} of $X.$ Then, since it follows from 
Lemma~\ref{pro'} that 
${Z\res(\eta\cup\tau)=(X\res(\eta\cup\tau))\,\cap\,
(Y\res(\eta'\cup\tau')\ares(\eta\cup\tau))},$ 
it suffices to verify that 
$x\res(\eta'\cup\tau')\in Y\res(\eta'\cup\tau').$ 
But $x\res\eta'\in Y\res\eta'$ by the choice of $x,$ and 
the same for $\tau',$ so that the required fact follows from 
property~\ref{indep} of $Y$.\qed

\ble
\label{gath'}
Let\/ $\i\in\bI,$ $X\in\ipe{\<\i},$ $X'\sq X$ is closed and 
nonempty, $X'\rsd{<\i}=X\rsd{<\i},$ and\/ $X'$ satisfies\/ 
{\rmt\ref{perf1}} and\/ {\rmt\ref{oz}} for this 
particular\/ $\i.$ Then\/ $X'\in\ipe{\<\i}$.
\ele
\proof Notice that requirement \ref{indep} is automatically 
satisfied in this case provided either $\xi$ or $\eta$ contains 
$\i$.\qed

\subsubsection*{Splitting of pre--conditions}   

We now demonstrate how a set in $\ipe\za$ can be 
splitted in a pair of smaller sets. Let $A\sq\cD$ be a set 
containing at least two different points. The largest finite 
sequence $r\in 2^{<\om}$ such that $r\subset a$ for all $a\in A$ 
is denoted by $\roo(A).$ We put $\ste(A)=\dom\roo(A)$ and define 
\dm
\spl(A,e)=\ans{a\in A:a(l)=e},\hspace{5mm}\hbox{where}
\hspace{3mm}l=\ste(A)\hspace{3mm}\hbox{and}\hspace{3mm} 
e=0\hspace{2mm}\hbox{or}\hspace{2mm}1\,.~\footnote
{\rmt\ Digits $0$ and $1$ will be denoted usually by letters 
$e$ and $d$.} 
\dm
Let now $X\in\ipe\za.$ Suppose that $\i\in\za.$ For any 
$y\in Y=X\rsd{<\i}$ the set $A(y)=D_{Xy}(\i)\sq\cD$ is perfect; 
therefore so are the sets $\spl(A(y),e),\,\,e=0,\,1.$ We define 
\dm
\spl(X,\i,e)=\ans{x\in X:x(\i)\in \spl(A(x\rsd{<\i}),e)}
\hspace{5mm}\hbox{for}\hspace{3mm}e=0,\,1\,.
\dm

\ble
\label{spl'}
Assume that\/ $\za\in\bI,$ $\i\in\za,$ and\/ $X\in\ipe\za$. 
Then\its
\ben
\item 
The sets $X_e=\spl(X,\i,e),\;e=0,\,1,$ belong to $\ipe\za$.\its

\item 
$X_0\res(\za\cap[\not\>\lsp\i])=X_1\res(\za\cap[\not\>\lsp\i])=
X\res(\za\cap[\not\>\lsp\i])$.\its

\item $X_0\rsd{\<\i}\cap X_1\rsd{\<\i}=\emptyset$.
\een
\ele
\proof First of all we note that since 
$X_e=X\cap(\spl(X\rsd{\<\i},\i,e)\ares\za),$ Assertion~\ref{less'} 
and Lemma~\ref{apro'} allow to consider only the case 
$\za=[\<\lsp\i].$ In this case items 2 and 3 become obvious, so 
we concentrate on item 1.

To prove that $X_e\in\ipe{\<\i},$ we put $Y=X\rsd{<\i},$ 
$D(y)=D_{Xy}(\i),$ and then define 
$Y_r=\ans{y\in Y:\roo(D(y))=r}$ for all 
$y\in Y$ and $r\in 2^{<\om}.$ It is implied by property {\ref{oz}} 
of $X$ that the sets $Y_r$ are clopen in $Y,$ and in fact there 
exist only finitely many nonempty sets $Y_r.$ Therefore the sets 
\dm
X_e={\textstyle\bigcup_r}\ans{x\in X:x\rsd{<\i}\in Y_r\cj 
x(\dom r)=e}\,,\hspace{5mm}e=0,\,1\,,
\dm
are clopen in $X,$ and $X_e\rsd{<\i}=X\rsd{<\i}.$ 
Furhermore condition \ref{perf1} for $X_e$ for the given $\i$ 
follows from the fact that nonempty intersections of perfect and 
clopen sets are perfect. Finally condition \ref{oz} for $X_e$ for 
the given $\i$ can be easily obtained from \ref{oz} for $X$ using 
the decomposition given by the last displayed 
formula.\qed
%
\newpage

\subsection{The forcing}
\label{forc}

The splitting procedure plays principal role in the complete 
definition of the notion of forcing. Let us start with several 
auxiliary definitions.

\bdf
Let $\za\sq\bI.$ A \dd\za{\it admissible\/} function is a function 
$\Phi:\om\;\lra\;\za$ taking each value $\i\in\za$ infinitely many 
times.\qed
\edf

\bdf
Assume that $\za\in\IS$ and $\Phi$ is a \dd\za admissible 
function. Let $X\in\ipe\za.$ We define a set $X[u]=X_\Phi[u]$ for 
all $u\in 2^{<\om}$ as follows.\its
\ben
\item $X[\La]=X$. ($\La$ is the empty sequence, the only member 
of $2^0$.)\its

\item If $m\in\om$ $u\in 2^m,$ and $X[u]$ has been defined, we 
put $X[u\we e]=\spl(X[u],\i,e),$ where $\i=\Phi(m),$ for 
$e=1,\,2$.\its
\een
For an infinite sequence $a\in 2^\om$ we define 
$X[a]=X_\Phi[a]=\bigcap_{m\in\om} X[a\res m]$. (Notice that 
$X[a]$ is nonempty by the compactness of $\can\za$.)\qed
\edf

\bcor
\label{itspl'}
If\/ $X\in\ipe\za$ then\/ 
$X_\Phi[u]\in\ipe\za$ for all\/ $u\in 2^{<\om}$.\qed
\ecor

\bdf
(Forcing conditions)\\[1pt]
Let $\za\in\IS.$ A set $X\in\ipe\za$ is {\it shrinkable\/} if for 
any \dd\za admissible function $\Phi$ and any $a\in 2^\om,$ the 
set $X_\Phi[a]$ contains only one point. 

We put $\pe\za=\ans{X\in\ipe\za:X\,\,\hbox{is shrinkable}}$ and 
$\pei=\pe\bI$.\qed
\edf

\noi
Obviously $X=\can\za$ is shrinkable (and belongs to $\pei$). 
We can easily prove that if $\za$ is wellfounded then every 
$X\in\ipe\za$ is shrinkable, so that $\pe\za=\ipe\za.$ On the 
other hand if $\za=\om^\ast$ (the order of negative integers) 
then the set
\dm
X=\ans{x\in\can\za:\forall\,\i\in\za\,(x(\i)(0)=0)\,\hbox{ or }\,
\forall\,\i\in\za\,(x(\i)(0)=1)}
\dm 
is not shrinkable: every set $X_\Phi[a]$ contains two different 
points. 

We now show that the lemmas already proved for pre--conditions 
remain valid for sets in $\pe\za$ --- {\it conditions\/} --- 
as well. 

\bass
\label{less}
If\/ $X\in\pe\za$ and\/ $\xi\in\IS,$ $\xi\sq \za,$ then\/ 
$Y=X\res \xi\in\pe\xi$.
\eass
\proof We have only to verify that $Y$ is shrinkable. Let $\Psi$ 
be a \dd\xi admissible function, and $b\in 2^\om.$ We want to 
prove that $Y_\Psi[b]$ is a singleton. Let $\Phi$ be a \dd\za 
admissible function such that $\Phi(2n)=\Psi(n)$ and 
$\Phi(2n+1)\in\za\setminus\xi$ for all $n.$ Let $a\in 2^\om$ be 
defined by $a(2n)=b(n)$ and $a(2n+1)=0$ for all $n.$ Then 
$X_\Phi[a]$ is a singleton since $X$ is shrinkable; on the other 
hand, $Y_\Psi[b]=X_\Phi[a]\res\xi$.\qed

\ble
\label{apro}
Suppose that\/ $\xi,\,\za\in\IS,$ $\xi\sq \za,$ $X\in\pe\za,$ 
$Y\in\pe\xi,$ and\/ $Y\sq X\res\xi.$ Then\/ $Z=X\cap(Y\ares\za)$ 
belongs to $\pe\za$. 
\ele
\proof To verify shrinkability let $\Phi$ be a \dd\za admissible 
function and ${a\in 2^\om}.$ To prove that $Z_\Phi[a]$ is a 
singleton, let ${\Om=\ans{k:\Phi(k)\in\xi}=\ans{o_m:m\in\om}}$ in 
the increasing order. We define $\Psi(m)=\Phi(o_m),$ so that 
$\Psi$ is \dd\xi admissible. We also put $b(m)=a(o_m).$ Then 
$Y_\Psi[b]$ is a singleton; on the other hand, $Y_\Psi[b]=
Z_\Phi[a]\res\xi$ because splittings $\spl(...,\i,e)$ with 
$\i\in\za\setminus\xi$ do not change projections on $\xi$. 

Thus at least $Z_\Phi[a]\res\xi=\ans{y}$ is a singleton. 

Let now $Y'=X\res\xi,$ so that $Y\sq Y'$ and $y\in Y'.$ By the 
shrinkability of $X,$ there exists some $a'\in 2^\om$ such that 
$y=x'\res\xi$ where $x'$ is the only element of $X_\Phi[a'].$ 
We now compose $c\in 2^\om$ from $a'$ and $a$ as follows: 
$c(k)=a'(k)$ for $k\in\Om$ and $c(k)=a(k)$ for 
$k\in\om\setminus\Om$.

Since $c\res\Om=a'\res\Om,$ we can easily prove that 
$X_\Phi[c]\res\xi=X_\Phi[a']\res\xi=\ans{y}=Z_\Phi[a].$ 
Furthermore since $c\res(\om\setminus\Om)=a\res(\om\setminus\Om),$ 
we obtain that in general $X_\Phi[c]=Z_\Phi[a],$ as required.\qed

\ble
\label{gath}
If\/ $X\in\pe{\<\i}$ in 
Lemma~\ref{gath'} then\/ $X'$ also 
belongs to $\pe{\<\i}$.
\ele
\proof $X'\rsd{<\i}=Y=X\rsd{<\i}$ is shrinkable by 
Lemma~\ref{less}. On the other hand, for any $y\in Y$ the set 
$D_{X'y}(\i)=\ans{x(\i):x\in X'}$ is converted to a singleton 
after infinitely many operations $\spl(...,\i,e)$\qed

\ble
\label{spl}
Assume that\/ $X\in\pe\za$ in Lemma~\ref{spl'}. Then the sets\/ 
$X_0$ and\/ $X_1$ also belong to $\pe\za$.
\ele
\proof The splitting of $X_e$ via an admissible $\Phi$ is equal to 
the splitting of $X$ itself via the function $\Psi$ defined so 
that $\Psi(0)=e$ and $\Psi(m+1)=\Phi(m)$.\qed

\bcor
\label{itspl}
If\/ $X\in\pe\za$ then\/ 
$X_\Phi[u]\in\pe\za$ for all\/ $u\in 2^{<\om}$.\qed
\ecor

We now come to the principal point which was perhaps the only 
reason for the introduction of shrinkability to the definition 
of forcing conditions. 

\bpro
\label{clop}
Assume that\/ $\za\in\IS,$ $X\in\pe\za,$ $X'\sq X$ is open in\/ 
$X,$ and\/ $x_0\in X'.$ There exists a clopen in\/ $X$ set\/ 
$X''\in\pe\za,$ $X''\sq X',$ containing $x_0$.
\epro
\proof Let $\Phi$ be a \dd\za admissible function. Then 
$\ans{x_0}=X_\Phi[a]$ for some (unique) $a\in 2^\om.$ By 
compactness there exists $m\in\om$ such that 
$X''=X_\Phi[a\res m]\sq X'$.\qed 

\newpage

\subsection{Finite iterations of splitting}
\label{ffin}

We shall exploit later the construction of sets in $\pei$ as 
$X=\bigcap_{m\in\om}\bigcup_{u\in 2^m} X_u,$ where every $X_u$ 
belongs to $\pei.$ Each level $\ang{X_u:u\in 2^m}$ of the given 
family of sets $X_u\in\pei$ should satisfy certain requirements 
which resemble the state as if we had defined the sets by 
iteration of splitting. This section introduces the requirements 
and presents some related lemmas. 

To specify the requirements which imply a good behaviour of the 
sets $X_u$ with respect to projections, we need to determine, for 
any pair of infinite binary sequences $u,\,v\in 2^{\<\om}$ of 
equal length, the largest initial segment $\xi=\xi[u,v]$ such that 
$X_u\res\xi$ should coincide with $X_v\res\xi$.

\bdf 
Let $\Phi$ be an \dd\bI admissible function. We define, for any 
pair of finite sequences $u,\,v\in 2^m,$ $m\in\om,$ an initial 
segment $\xi[u,v]=\xi_\Phi[u,v]\in \IS$ by induction as follows.\its
\ben
\itla{xi0} 
$\xi[\La,\La]=\bI$. ($\La$ is the empty sequence.)\its

\itla{xieq} 
$\xi[u\we d,v\we e]=\xi[u,v]$ provided $d=e\in\ans{0,1}$.\its

\itla{xiOm} 
Assume that $u$ and $v$ belong to $2^n,$ 
$d,\,e\in\ans{0,1}$ are different, and $\Phi(n)=\i\in\bI.$
Then $\xi[u\we d,v\we e]=\xi[u,v]\cap[\not\>\lsp\i]$.\its\qed 
\een
\edf

\bdf
Let $\Phi$ be an \dd\bI admissible function. 
A \dd\Phi{\it splitting system\/} of order $m$ is an indexed 
family $\ang{X_u:u\in 2^{\<m}}$ of sets $X_u\in \pei$ such that\its
\ben
\def\theenumi{\rmt{S-}\arabic{enumi}}
\itla{spli}
$X_{u\we e}\sq \spl(X_u,\Phi(n),e)$ whenever 
$u\in 2^n,$ $n<m,$ and $e\in\ans{0,1}$.\its

\itla{prct} 
$X_u\res\xi[u,v]=X_v\res\xi[u,v]$ for any pair of 
$u,\,v$ of equal length $\<m$.\its

\itla{aprct} 
If $\i\in\bI,$ $\i\not\in\xi[u,v],$ then 
$X_u\rsd{\<\i} \cap X_v\rsd{\<\i}=\emptyset$.\qed
\een
\edf 

\noi In particular it easily follows from lemmas \ref{spl'} and 
\ref{spl} that for all $m$ and $X\in\pei$ the system defined by 
$X_u=X_\Phi[u]$ for each $u\in 2^{\<m},$ is a \dd\Phi splitting 
system.

We consider two ways how an existing splitting system can be 
transformed to another splitting system. One of them treats 
the case when we have to change one of sets to a smaller set, 
the other one is an expansion to the next level.

It is assumed that a \dd\bI admissible function $\Phi$ is fixed 
and $\xi[u,v]=\xi_\Phi[u,v]$.

\ble
\label{suz}
Assume that\/ $\ang{X_u:u\in 2^{\<m}}$ is a splitting system, 
$u_0\in 2^m,$ and\/ $X\in\pei,$ $X\sq X_{u_0}.$ We re--define\/ 
$X_u$ by\/ 
$X'_u=X_u\cap (X\res\xi[u,u_0]\ares\bI)$ for all\/ $u\in 2^m.$  
Then the re--defined~\footnote
{\rmt\ Notice that $X'_{u_0}=X$.} 
family is again a splitting system. 
\ele
\proof All sets $X'_u$ belong to $\pei$ by Lemma~\ref{apro}. 
Therefore we have to check only condition~\ref{prct}. Thus let 
$u,\,v\in 2^m,$ $\xi=\xi[u,v].$ We prove that 
$X'_u\res\xi=X'_v\res\xi$. Let in addition $\xi_u=\xi[u,u_0]$ 
and $\xi_v=\xi[v,u_0].$ Then 
\dm
X'_u\res\xi=X_u\res\xi\,\cap\,(X_0\res(\xi\cap\xi_u)\ares\xi)
\hspace{4mm}\hbox{and}\hspace{4mm}
X'_v\res\xi=X_v\res\xi\,\cap\,(X_0\res(\xi\cap\xi_v)\ares\xi)
\dm
by Lemma~\ref{pro'}. Thus it remains to prove that 
$\xi\cap\xi_u=\xi\cap\xi_v.$ Assume that on the contrary 
$\i\in\xi\cap\xi_u$ but $\i\not\in\xi_v.$ The 1st assumption 
implies $X_{u_0}\rsd{\<\i}=X_u\rsd{\<\i}=X_v\rsd{\<\i}$ by 
condition~\ref{prct}, but the 2nd one implies that 
$X_{u_0}\rsd{\<\i}\cap X_v\rsd{\<\i}=\emptyset,$ 
contradiction.\qed

\ble
\label{pand}
Every splitting system\/ $\ang{X_u:u\in 2^{\<m}}$ can be expanded 
to the next level by adjoining appropriate sets\/ $X_{u'},$ 
$u'\in 2^{m+1}$. 
\ele
\proof Let $\Phi(m)=\i\in\bI.$ We define 
$X_{u\we e}=\spl(X_u,\i,e)$ for all $u\in 2^m$ and $e=0,\,1.$ It 
suffices to prove conditions \ref{prct} and \ref{aprct}. 
Let $u'=u\we d$ and $v'=v\we e$ belong to $2^{m+1}$. We put 
$Y=X_u\res\xi=X_v\res\xi,$ where $\xi=\xi[u,v]$.\vom

{\it Case 1\/}: $\i\not\in\xi.$ Then by definition 
$\xi=\xi[u',v']$ as well. Lemma~\ref{spl'}  
immediately gives $X_{u'}\res\xi=Y=X_{v'}\res\xi.$ This proves 
\ref{prct}. On the other hand, if $\i\not \in \xi$ then already 
$X_u\rsd{\<\i}\cap X_v\rsd{\<\i}=\emptyset,$ and we have 
\ref{aprct}.\vom

{\it Case 2\/}: $\i\in\xi$ and $d=e,$ say $d=e=0.$ Then again 
$\xi=\xi[u\we d,v\we e]=\xi[u,v].$ We obtain
$X_{u'}\res\xi=\spl(Y,\i,0)=X_{v'}\res\xi,$ as required.\vom

{\it Case 3\/}. $\i\in\xi$ and $d\not=e,$ say $d=0,$ $e=1.$ 
Then $\xi'=\xi[u\we d,v\we e]=\xi\cap[\not\>\lsp\i],$ and we 
obtain 
$X_{u'}\res\xi'=\spl(Y,\i,0)\res\xi'=\spl(Y,\i,1)\res\xi'=
X_{v'}\res\xi',$ as required. (The middle equality follows from 
Lemma~\ref{spl'} for $\za=\xi.$) To check~\ref{aprct} for some 
$\j,$ note 
that if $\j\not\in\xi'$ then either already $\j\not\in\xi$ or 
$\j\>\i.$ In the latter case we use Lemma~\ref{spl'} again 
to obtain $\spl(Y,\i,0)\rsd{\<\i}\cap \spl(Y,\i,1)\rsd{\<\i}=
\emptyset.$ But $X_{u'}\rsd{\<\i}=\spl(Y,\i,0)\rsd{\<\i}$ and 
$X_{v'}\rsd{\<\i}=\spl(Y,\i,1)\rsd{\<\i}$ because $\i\in\xi$ and 
$Y=X_u\res\xi=X_v\res\xi$.\qed

\newpage
\subsection{The forcing is homogeneous}
\label{hom}

The following theotem shows that the forcing $\pe\za$ is quite 
homogeneous. This fact will be used in the proof of Theorem~\ref{ge} 
in Section~\ref{seqge} but not earlier.
 
\bte
\label{thom}
Suppose that\/ $\za\in\IS.$ Let\/ $X,\,Y\in\pe\za.$ There exists 
a homeomorphism\/ $F:X\,\lra\,Y$ such that\its
\ben
\itla{h1} For all\/ $\xi\in\IS,$ $\xi\sq\za,$ and $x,\,y\in X,$ 
$x\res\xi = y\res\xi$ iff\/ $F(x)\res\xi = F(y)\res\xi$.\its

\itla{h2} If\/ $\xi\in\IS,$ $\xi\sq\za,$ and $X\res\xi=Y\res\xi,$ 
then\/ $x\res\xi=F(x)\res\xi$ for all\/ $x\in X$.\its

\itla{h3} For each\/ $X'\sq X$ 
we have\/ $X'\in\pe\za$ iff\/ $Y'=F\ima X'\in\pe\za$.\\
{\rmt ($F\ima X=\ans{F(x):x\in X}$ is the image of $X$ via $F$.)}
\een
\ete
\proof Let us fix a \dd\za admissible function $\Phi.$ We recall 
that the initial segment $\xi[u,v]=\xi_\Phi[u,v]\in\IS$ is defined 
for $u,\,v\in 2^n$ by induction on $n$ in Section~\ref{ffin}. We 
set $\xi[a,b]=\bigcap_n\xi[a\res n,b\res n]$ for $a,\,b\in 2^\om;$ 
then $\xi[a,b]\in\IS$ as well. 

We defined sets $X[u]=X_\Phi[u]\;\;(u\in 2^{<\om})$ in the end of 
Section~\ref{splt}; all of them belong to $\pe\za$ by 
Corollary~\ref{itspl}. Lemma~\ref{pand} implies:
\ben
\def\theenumi{{\hskip1pt\rmt{s-}\arabic{enumi}\hskip1pt}}
\itla{s1}\msur
$X[u\we e]\sq\spl(X[u],\Phi(n),e)$ whenever 
$u\in 2^n,$ $n<m,$ and $e\in\ans{0,1}$.\its

\itla{s2}\msur
$X[u]\res\xi[u,v]=X[v]\res\xi[u,v]$ for any pair of 
$u,\,v$ of equal length $m\in\om$.\its

\itla{s3} 
If $\i\in\bI,$ $\i\not\in\xi[u,v],$ then 
$X[u]\rsd{\<\i}\cap X[v]\rsd{\<\i}=\emptyset$.\qed
\een
Since $X,\,Y\in\pe\za,$ the sets $X[a]=\bigcap_nX[a\res n]$ and 
$Y[a]=\bigcap_nY[a\res n]$ contain one point each, resp. $x_a$ 
and $y_a,$ for every $a\in 2^\om.$ It follows from \ref{s2} 
and \ref{s3} that\its
\bit
\item[$(\ast)$] $\msur x_a\res \xi[a,b]=x_b\res \xi[a,b]$ \ 
for all \ $a,\,b\in 2^\om.$ \ \ If \ $\i\not\in\xi[a,b]$ \ then \ 
$x_a\rsd{\<\i}\not=x_b\rsd{\<\i}$.\its
\eit
It follows that $x_a\res\xi= x_b\res\xi$ iff $\xi\sq\xi[a,b]$ for 
each pair of $a,\,b\in 2^\om$ and every $\xi\in\IS,\;\,\xi\sq\za;$ 
the same is true for $y_a,\,y_b.$  
Taking $\xi=\za,$ we define a homeomorphism $F:X\,$onto$\,Y$ by 
$F(x_a)=y_a$ for all $a\in 2^\om.$ For each 
$\xi\in\IS,\;\,\xi\sq\za$ we obtain an associated homeomorphism 
$F_\xi:X\res\xi\,$onto$\,Y\res\xi$ which satisfies\its
\bit
\item[$(\dag)$] 
$F_\xi(x\res\xi)=F(x)\res\xi$ \ for all \ $x\in X$.\its
\eit
Let us prove that $F$ is as required. Item~\ref{h1} 
of the theorem immediately follows from $(\dag).$ To verify 
item~\ref{h2}, one easily proves that $X[u]\res\xi=Y[u]\res\xi$ for 
all $u\in 2^{<\om}$ by induction on the length of $u,$ provided 
$X\res\xi=Y\res\xi.$ It remains to check item~\ref{h3}.

Let $X'\sq X,\;\,X'\in\pe\za;$ we have to prove that 
$Y'=F\ima X'\in\pe\za$ as well.\vspace{1mm}

{\it We check requirement\/ \ref{perf1}\/}. 
Thus let $\i\in\za$ and $y'\in Y'\rsd{<\i}.$ To prove that 
$\hat Y=D_{Y'y'}(\i)$ is a perfect set, let 
$y'=F_{<\i}(x');\;x'\in X'\rsd{<\i}.$ Then the set 
$\hat X=D_{X'x'}(\i)$ is perfect. Let $\hat x\in\hat X;$ then 
$\hat x=x(\i)$ for some $x\in X'$ such that $x\rsd{<\i}=x'.$ 
One can define $y=F(x)$ and $\hat y=y(\i)$ -- then 
$y\in Y',$ $y\rsd{<\i}=y',$ and $\hat y\in \hat Y$ by $(\dag).$ 
It also follows from $(\dag)$ that in fact $\hat y$ depends 
only on $\hat x$ but not on the entire $x;$ let 
$\hat y=\hat F(\hat x).$ We conclude that $\hat F$ is a 
homeomerphism $\hat X$ onto $\hat Y,$ as required.\vspace{1mm}

{\it We check requirement\/ \ref{oz}\/}. Suppose that 
$\i\in\za$ and $Y^\ast\sq Y'$ is clopen in $Y';$ we have to 
prove that the set $Y^\ast\rsd{<\i}$ is open in $Y'\rsd{<\i}.$ 
(This is more than \ref{oz} asserts, but here this is more 
convenient.) 
We put $X^\ast=F^{-1}(Y^\ast);$ then $X^\ast$ is clopen in 
$X'$ because $F$ is a homeomorphism, and 
$Y^\ast\rsd{<\i}=F_{<\i}\ima (X^\ast\rsd{<\i})$ by $(\dag);$ 
therefore it suffices to prove that 
$X^\ast\rsd{<\i}$ is clopen in $X'\rsd{<\i}$. 

Since each $X'[a]=X'_\Phi[a]\;\;(a\in 2^\om)$ is a singleton by 
the shrinkability, there exists a number $n\in\om$ such that we 
have $X'[u]\sq X^\ast$ or $X'[u]\cap X^\ast=\emptyset$ for every 
$u\in 2^n.$ It follows from \ref{s2} and \ref{s3} above that 
the projections $X'[u]\rsd{<\i}$ and $X'[v]\rsd{<\i}$ are either 
equal or disjoint for each pair of $u,\,v\in 2^n.$ This easily 
implies the required fact.\vspace{1mm}

{\it We check requirement\/ \ref{indep}\/}. Let 
$\xi,\,\eta\in\IS,$ $\xi\cup\eta\sq\za,$ $y'\in Y'\res\xi,$ 
$y''\in Y'\res \eta,$ and 
$y'\res (\xi\cap \eta)=y''\res (\xi\cap \eta);$ we have to prove 
that $y'\cup y''\in Y'\res (\xi\cup\eta).$ Using $(\dag),$ we 
see that $y'=F_\xi(x')$ and $y''=F_\eta(x''),$ where 
$x'\in X'\res\xi,$ $x''\in X'\res \eta,$ and 
$x'\res (\xi\cap \eta)=x''\res (\xi\cap \eta).$ Then 
$x'\cup x''\in X'\res (\xi\cup\eta)$ because $X'\in\pe\za.$ We 
conclude that 
$y'\cup y''=
F_{\xi\cup\eta}(x'\cup x'')\in Y'\res (\xi\cup\eta),$ 
as required.\vspace{1mm}

{\it We prove that $Y'$ is shrinkable\/}. Since $X'$ is shrinkable, 
it suffices to verify that $Y'[u]=F\ima X'[u]$ for all 
$u\in2^{<\om}.$ It follows from $(\dag)$ that the problem can be 
reduced to one--dimentional setting. 

Let $P$ be a perfect subset of $\cD=2^\om.$ We set $P[\La]=P$ and 
$P[u\we e]=\spl(P[u],e)$ for all $u\in 2^{<\om}$ and $e\in\ans{0,1}$ 
(see the end of Section~\ref{splt}); thus a perfect set $P[u]$ is 
defined for 
all $u\in 2^{<\om}.$ Obviously $P[a]=\bigcap_n P[a\res n]$ is a 
singleton, say $p_a,$ for every $a\in 2^\om.$ Thus 
$a\,\longmapsto\,p_a$ is a homeomorphism $2^\om$ onto $P$. 

Suppose that $P,\,Q$ is a pair of perfect subsets of $\cD,$ 
with the associated homeomorphisms $a\longmapsto p_a$ and 
$a\longmapsto q_a.$ One 
defines a special homeomorphism $f=f_{PQ}:P\,$onto$\,Q$ by 
$f(p_a)=q_a$ for all $a\in 2^\om.$ Now, the abovementioned 
one--dimentional assertion 
in the verification of shrinkability, is as follows:\its
\bit
\item[$(\ddag)$] {\it prove that\/ $f\ima (P'[u])=(f\ima P')[u]$ for 
every perfect\/ $P'\sq P$ and all\/ $u\in2^{<\om}$.}\its
\eit
We prove $f\ima (P'[u])=(f\ima P')[u]$ for a fixed perfect $P'\sq P$ 
by induction on the length of $u\in 2^{<\om}.$ The case $u=\La$ is 
clear. We now suppose that $f\ima (P'[u])=Q'[u],$ where 
$Q'=f\ima P',$ and prove 
$f\ima (P'[u\we e])=Q'[u\we e]\,,$ $e=0,\,1$.

Let $s\in 2^{<\om}$ be the maximal sequence such that 
$P'[u]\sq P[s].$ By definition $f=f_{PQ}$ maps $P[s]$ onto 
$Q[s];$ by the assumption $f$ also maps $P'[u]$ onto 
$Q'[u].$ We observe that $Q'[u]\sq Q[s]\,;$ moreover, 
$Q'[u]\not\sq Q[s\we e]$ for any $e=0,\,1$ (otherwise we would 
get $P'[u]\sq P[s\we e]\,,$ contradiction with the choice of $s$). 
It follows that
\dm
P'[u\we e]=P'[u]\cap P[s\we e]
\hspace{8mm}\hbox{and}\hspace{8mm}
Q'[u\we e]=Q'[u]\cap Q[s\we e]
\dm
for $e=0,\,1.$ 
This implies $f\ima (P'[u\we e])=Q'[u\we e]\,,$ as required.\qed

\newpage

\subsection{Fusion technique}
\label{finf}

We prove here a version of {\it fusion lemma\/} which fits to this 
form of iterated Sacks forcing. 

An \dd\bI admissible function $\Phi$ is fixed; all notions of the 
preceding section are treated in the sense of this function $\Phi$. 

\bdf
An indexed family of sets $X_u\in\pei,$ $u\in 2^{<\om},$ 
is a {\it fusion sequence\/} if first, for every $m\in \om$ 
the subfamily $\ang{X_u:u\in 2^{\<m}}$ is a splitting system, 
and second,\its
\ben
\def\theenumi{\rmt{S-}\arabic{enumi}}
\setcounter{enumi}{3}
\itla{di} 
For any $\vep>0$ there exists $m\in\om$ such that $\diam X_u<\vep$ 
for all $u\in 2^m.$  
(A Polish metric on $\cani$ is fixed.)\qed
\een
\edf

\bte
\label{fut}
Let\/ $\ang{X_u:u\in 2^{<\om}}$ be a fusion sequence. Then\/ 
$X=\bigcap_{m\in\om}\bigcup_{u\in 2^m} X_u$ belongs to 
$\pei$.
\ete
\proof The intersection $\bigcap_m X_{a\res m}$ contains a single 
point $x_a\in X$ for any $a\in 2^\om=\cD$ by \ref{di}, and the 
map $a\longmapsto x_a$ is continuous. Let us define 
$\xi[a,b]=\bigcap_{m\in\om}\xi[a\res m,b\res m].$ In particular 
$\xi[a,b]=\bI$ iff $a=b.$ It follows from conditions \ref{prct} and 
\ref{aprct} that\its
\bit
\item[$(\ast)$] $\msur x_a\res \xi[a,b]=x_b\res \xi[a,b]$ 
for all $a,\,b\in 2^\om.$ If $\i\not\in\xi[a,b]$ then 
$x_a\rsd{\<\i}\not=x_b\rsd{\<\i}.$\its
\eit

\noi We check condition~\ref{perf1}. Let $\i\in \bI$ and 
$y\in X\rsd{<\i},$ $p\in D=D_{Xy}(\i).$ For a fixed $k\in\om,$ we 
shall find a point $q\in D,$ $q\not=p,$ which satisfies 
$q\res k=p\res k.$ Let $x\in X$ be such that $p=x(\i)$ and 
$x\rsd{<\i}=y.$ Let $x=x_a;\;a\in 2^\om.$ First of all we fix 
$m\in \om$ such that $x'(\i)\res k=x''(\i)\res k$ for all 
$u\in 2^m$ and $x',\,x''\in X_u.$ Let $b\in 2^\om$ be defined by 
$b(n)=a(n)$ for all $n\in\om$ with the exception of the first 
$n>m$ such that $\Phi(n)=\i,$ where we set $b(n)=1-a(n).$ Let 
$q=x_b(\i)$. Then\its
\bit
\item[$1)$] 
taking $u=a\res m=b\res m,$ we get $p\res k= q\res k$;\its

\item[$2)$] evidently $[<\lsp\i]\sq\xi[a,b],$ therefore 
$x_b\rsd{<\i}=x_a\rsd{<\i}=y,$ so that $q\in D$;\its

\item[$3)$]\msur 
$q=x_b(\i)\not=x_a(\i)=p$ because $a(n)\not=b(n)$ for some 
$n$ such that $\Phi(n)=\i$.\its
\eit
(In the last two items, we use $(\ast)$.)\vom

We check condition~\ref{oz} for the set $X.$ Let $\i\in \bI$ and 
$G\sq\cD$ be an open set. We define $X'=\ans{x\in X:x(\i)\in G}$ 
and prove that $X'\rsd{<\i}$ is open in $X\rsd{<\i}$. 

We can assume that $G$ is in fact clopen. Then by \ref{di} 
there exists $m\in\om$ such that for every $u\in 2^m$ either 
$X_u(\i)\sq G$ or $X_u(\i)\cap G=\emptyset.$ Let 
$U=\ans{u\in 2^m:X_u(\i)\sq G}.$ Notice that in accordance with 
conditions \ref{prct} and \ref{aprct}, for any pair 
$u,\,v\in2^m,$ either $X_u\rsd{<\i}=X_v\rsd{<\i}$ or 
$X_u\rsd{<\i}\cap X_v\rsd{<\i}=\emptyset.$ Let $V$ be the set of 
all $v\in 2^m$ such that $X_u\rsd{<\i}\cap X_v\rsd{<\i}=\emptyset$ 
for all $u\in U.$ There exists a clopen set $C\sq\can{<\i}$ which 
separates $A=\bigcup_{u\in U}X_u\rsd{<\i}$ from 
$B=\bigcup_{v\in V}X_v\rsd{<\i}.$ It remains to verify that 
$X'\rsd{<\i}=X\rsd{<\i}\cap C$.

Let $x\in X',$ that is, $x=x_a\in X$ and $x(\i)\in G.$ Then $u=
a\res m\in U,$ therefore $x\rsd{<\i}\in A$ and $\in C.$ Let for 
the converse $x=x_a\in X$ and $x\rsd{<\i}\in C;$ we have to find 
$x'\in X'$ such that $x'\rsd{<\i}=x\rsd{<\i}.$ Notice that $v=
a\res m$ does not belong to $V$ (although may not belong to $U$ 
either). Therefore $X_v\rsd{<\i}$ is equal to $X_u\rsd{<\i}$ for 
some $u\in U.$ Let $b\in 2^\om$ be defined by $b\res m=u$ and 
$b(n)=a(n)$ for $n\>m.$ Then, since $X_v\rsd{<\i}=X_u\rsd{<\i},$ 
we obtain $X_{a\res n}\rsd{<\i}=X_{b\res n}\rsd{<\i}$ for all 
$n>m.$ Therefore, $x'=x_b$ satisfies $x'\rsd{<\i}=x\rsd{<\i}.$ On 
the other hand, $x'\in X'$ because $b\res m=u\in U$.\vom

We check \ref{indep}. Let $\xi,\,\eta\in\IS,$ and points 
$x'=x_{a'},\,x''=x_{a''}$ of $X$ be such that 
$x'\res(\xi\cap\eta)=x''\res(\xi\cap\eta);$ then in particular 
$\xi\cap\eta\sq\xi[a',a'']$ by $(\ast)$. We have to find 
$x=x_a\in X$ satisfying $x\res\xi=x'\res\xi$ and 
$x\res\eta=x''\res\eta$.  

To obtain the required $a\in 2^\om,$ we define the values 
$a(m)\in\ans{0,1}$ using induction on $m.$ Assume that we have 
defined $a\res m,$ and define $a(m)$. Let $\Phi(m)=\i\in\bI$.\vom

{\it Case 1\/}: $\i\not\in\xi\cup\eta.$ We define $a(m)$ 
arbitrarily, say $a(m)=0$ in this case.

{\it Case 2\/}: $\i\in\xi\setminus\eta$ or 
$\i\in\eta\setminus\xi.$ We put resp. $a(m)=a'(m)$ or 
$a(m)=a''(m)$.

{\it Case 3\/}: $\i\in\xi\cap\eta.$ Then $a'(m)=a''(m)$ since 
otherwise we would have $\i\not\in\xi[a',a''].$ We put 
$a(m)=a'(m)=a''(m)$.\vom

This definition implies $\xi\sq\xi[a,a']$ and 
$\eta\sq\xi[a,a''].$ Then $x_a\res\xi=x_{a'}\res\xi$ and 
$x_a\res\eta=x_{a''}\res\eta$ by $(\ast),$ as required.

Finally to check shrinkability we note that $X_\Phi[u]=
X\cap X_u$ (can be easily proved by induction), so that the 
property follows from~\ref{di}.\qed

\bcor
\label{prebor}
Let\/ $\xi\in\IS,$ $X\in\pe\xi,$ and\/ $C_m\sq\can\xi$ be closed 
for each $m\in\om.$ There exists\/ $Y\in\pe\xi,$ $Y\sq X$ such 
that\/ $C_m\cap Y$ is clopen in\/ $Y$ for every\/ $m$.
\ecor
\proof We can assume that $\xi=\bI$ (if not replace $C_m$ by 
$C_m\ares\bI$). It follows from Proposition~\ref{clop} that for 
any $m$ and any $X'\in\pei$ there exists $Y'\in\pei,$ $Y'\sq X',$ 
such that either $Y'\sq C_m$ or $Y'\cap C_m=\emptyset.$ Therefore 
we can define, using lemmas \ref{suz} and \ref{pand}, a fusion 
sequence $\ang{X_u:u\in 2^{<\om}}$ of sets $X_u\in\pei$ such that 
$X_\La=X$ and either $X_u\sq C_m$ or $X_u\cap C_m=\emptyset$ 
whenever 
$u\in 2^m.$ Let $Y=\bigcap_{m\in\om}\bigcup_{u\in 2^m}X_u$.\qed

\bcor
\label{bor}
Assume that\/ $\xi\in\IS,$ $X\in\pe\xi,$ and\/ $B\sq\can\xi$ 
is a set of a finite Borel level. There exists\/ $Y\in\pe\xi,$ 
$Y\sq X$ such that either\/ $Y\sq B$ or\/ $Y\cap B=\emptyset$. 
\ecor
\proof\footnote
{\rmt\ In fact this is true for all Borel sets $B$ but needs 
more elaborate reasoning.} \ 
Let $B$ be defined by a finite level Borel scheme (countable 
unions plus countable intersections) from closed sets $
C_m,\;m\in\om.$ The preceding corollary shows that there exists 
$X'\in\pe\xi,$ $X'\sq X$ such that every $X'\cap C_m$ is clopen 
in $X'.$ Thus the Borel level can be reduced.\qed\vspace{3mm}

The results already obtained are in fact sufficient to prove the 
first part of Theorem~\ref{m}. However to handle the degrees of 
constructibility in \dd\pei generic extensions we need to conduct a 
more detailed analysis concentrated on continuous functions. 

\newpage%
\subsection{Reducibility of continuous functions}
\label{cont}

This section provides analysis of the behaviour of continuous 
functions defined on sets in $\pei$ from the point of view of 
certain reducibility. 

\bdf
For $\xi\sq\bI,$ $\cnt\xi$ is the set of all continuous 
functions $F:\can\xi\;\lra\;\om^\om.$ We put $\cnti=\cnt\bI.$ 
Let $F\in\cnti$.\its
\ben
\itla{redu}\msur 
$F$ is {\it reducible\/} to 
$\xi\in\IS$ on a set $X\sq\cani$ if for all $x,\,y\in X$ such that 
$x\res\xi=y\res\xi$ we have $F(x)=F(y)$.\its

\itla{capt}\msur 
$F$ {\it captures\/} $\i\in\bI$ on $X$ 
if for all $x,\,y\in X$ such that $F(x)=F(y)$ we have 
$x(\i)=y(\i)$.\qed
\een
\edf

\brem
\label{redr}
It follows from the compactness of the spaces we consider, that 
if $X$ is closed then in item~\ref{redu} there exists a function 
$H\in\cnt\xi$ such that $F(x)=H(x\res\xi)$ for all $x\in X,$ while 
in item~\ref{capt} there exists a continuous function 
$E:\cD\;\lra\;\cD$ such that $x(\i)=E(F(x))$ for all $x\in X$.\qed
\erem

Theorem~\ref{m} contains four items, \ref{m1+} through \ref{m3}, 
concerning \dd\gM constructibility of reals in the extension. We 
shall obtain those assertions as corollaries of the following 
theorem.

\bte
\label{mm}
Assume that\/ $\xi\in\IS,$ $X\in\pei,$ $F\in\cnti.$ Then\its
\ben
\itla{mm1+} 
If\/ $\i,\,\j\in\bI$ and\/ $\i<\j$ then  
there exists\/ $X'\in\pei,\;\;X'\sq X,$ such that the 
co-ordinate function $C_\j$ defined on\/ $\cani$ by 
$C_\j(x)=x(\j)$ captures\/ $\i$ on $X'$.\its

\itla{mm2} 
If\/  $\i\in\bI\setminus\xi$ and\/ $F$ is reducible to\/ 
$\xi$ on\/ $X$ then\/ $F$ does not capture\/ $\i$ on\/ 
$X$.\its

\itla{mm2+} 
Suppose that\/  
$F$ satisfies the property that for all\/ $X'\in\pei,$ 
$X'\sq X,$ and\/ $\i\in\xi$ there exists\/ $X''\in\pei,$ 
$X''\sq X'$ such that\/ $F$ captures\/ $\i$ on\/ $X''.$ Then there 
exists\/ $Y\in\pei,$ $Y\sq X$ such that\/ $F$ captures each\/ 
$\i\in\xi$ on\/ $Y$.\its

\itla{mm3} 
There exists\/ $X'\in\pei,$ $X'\sq X$ satisfying exactly one from 
the following two assertions$:$\hspace{4mm} 
\begin{minipage}[t]{0.6\textwidth}
{\rmt(a)} $F$ is reducible to\/ $\xi$ on\/ $X'$, \hfill  
or\vspace{2mm}

{\rmt(b)} $F$ captures some\/ $\i\in\bI\setminus\xi$ on\/ $X'$.
\end{minipage}
\hfill $\,$
\een
\ete
\proof We begin from a few technical lemmas, then come 
to the theorem.

\ble
\label{inter}
If\/ $F$ is reducible to both\/ $\xi$ and\/ $\eta$ on\/ 
$X\in\pei$ then\/ $F$ is reducible to\/ $\za=\xi\cap\eta$ on 
$X$.
\ele
\proof Let, on the contrary, $x,\,y\in Y$ satisfy $x\res\za=
y\res\za$ but $F(x)\not=F(y).$ Then by property~\ref{indep} 
of $X$ there exists $z\in X$ such that $z\res\xi=x\res\xi$ and 
$z\res\eta=y\res\eta.$ We obtain $F(x)=F(z)=F(y),$ 
contradiction.\qed

\ble
\label{l1}
Assume that\/ $\xi\in\IS,$ $\i\in\bI\setminus\xi,$ the sets 
$X_1$ and\/ $X_2$ belong to $\pei,$ and\/ $X_1\res\xi=X_2\res\xi.$ 
Then either $F$ is reducible to $\xi$ on $X_1\cup X_2$ or there 
exist sets\/ $X_1',\,X_2'\in\pei,$ $X_1'\sq X_1$ and\/
$X_2'\sq X_2,$ such that again\/ $X'_1\res\xi=X'_2\res\xi,$ and\/ 
$F\ima X'_1\cap F\ima X'_2=\emptyset$.~\footnote
{\rmt\ We recall that $F\ima X$ is the image of $X$ via $F$.}
\ele
\proof We suppose that $F$ is not reducible to $\xi$ on 
$X_1\cup X_2$ and prove the ``or'' alternative. By the 
non--reducibility assumption there exist points 
$x_1,x_2\in X_1\cup X_2$ such that $x_1\res\xi= x_2\res\xi$ and 
$F(x_1)\not=F(x_2).$ By the continuity of $F$ there exist 
clopen neighbourhoods $U_1$ and $U_2$ of $x_1$ and $x_2$ such that 
$F\ima U_1\cap F\ima U_2=\emptyset.$ Proposition~\ref{clop} gives  
a set $X''_1\in\pei,$ $X''_1\sq X_1\cap U_1,$ containing $x_1$. 

By Lemma~\ref{apro} the set 
$X''_2=X_2\cap (X''_1\res\xi\ares\bI)$ belongs to $\pei,$ and 
contains $x_2$ since $x_1\res\xi=x_2\res\xi.$ By 
Proposition~\ref{clop} again, there exists a set $Z_2\in\pei,$ 
$Z_2\sq X''_2\cap U_2.$ Now putting 
$Z_1=X''_1\cap (Z_2\res\xi\ares\bI)$ we get sets 
$Z_1,\,Z_2\in\pei$ such that $Z_1\sq X_1,$ $Z_2\sq X_2,$ 
$Z_1\res\xi=Z_2\res\xi$ and 
$F\ima Z_1\cap F\ima Z_2=\emptyset$.\qed

\ble
\label{l2}
If\/ $F\in\cnti,$ $X\in\pei,$ and\/ $\i\in\bI,$ then there 
exists\/ $X'\in\pei,$ $X'\sq X$ such that either\/ $F$ is 
reducible to\/ $[\not\>\lsp\i]$ on\/ $X',$ or\/ $F$ 
captures\/ $\i$ on $X'$.
\ele
\proof Let us assume that a set $X'$ of the ``either'' type does 
not exist. To prove the existence of $X'$ of the ``or'' type, 
we fix an \dd\bI admissible function $\Phi$ and put 
$\xi[u,v]=\xi_\Phi[u,v]$ for every pair of finite sequences 
$u,\,v\in 2^{<\om}$ of equal length. The notions of splitting 
system and fusion sequence are understood in the sense of $\Phi$.

Using Lemma~\ref{l1} and Proposition~\ref{clop}, we define a 
fusion sequence $\ang{X_u:u\in 2^{<\om}}$ with $X_\La=X,$ 
satisfying the following condition:\its
\bit
\item [$(\star)$] If $m\in \om$ and $u,\,v\in 2^m$ then either 
$(1)$ $F$ is reducible to $\xi[u,v]$ on the set $X_u\cup X_v$, \ or 
$(2)$ $F\ima X_u\cap F\ima X_v=\emptyset$.
\its
\eit
Indeed we put $X_\La=X,$ as indicated. Assume that sets $X_u,$ 
$u\in 2^m,$ have been defined, and $\i=\Phi(m).$ We first set 
$Y_{u\we e}=\spl(X_u,\i,e)$ for all $u\in 2^m$ and $e=0,\,1,$ 
obtaining a splitting system $\ang{Y_{u'}:u'\in 2^{m+1}}$ (see 
Lemma~\ref{pand}). At the next step we consider consequtively 
all pairs $u',v'\in 2^{m+1}$ and reduce sets $Y_{u'}$ and 
$Y_{v'}$ using Lemma~\ref{l1} for $\xi=\xi[u,v].$ Let $X'_{u'}$ 
and $X'_{v'}$ be the pair of sets we obtain; in particular 
$X'_{u'}\res\xi=X'_{v'}\res\xi,$ and either $F$ is reducible 
to $\xi$ on $X'_{u'}\cup X'_{v'}$ or 
$F\ima X'_{u'}\cap F\ima X'_{v'}=\emptyset$.

Then we set 
$Z_{w'}=X'_{w'}\cap (X'_{u'}\res\xi[w',u']\ares\bI),$ so that 
$\ang{Z_{w'}:w'\in 2^{m+1}}$ is a splitting system by 
Lemma~\ref{suz}. It is essential that since 
$X'_{u'}\res\xi=X'_{v'}\res\xi,$ we have $X'_{v'}\sq Z_{v'}.$ 
This allows to repeat the reduction: let 
$Y'_{w'}=Z_{w'}\cap (X'_{v'}\res\xi[w',u']\ares\bI),$ which 
gives again a splitting system of sets such that 
$Y'_{u'}=X'_{u'}$ and $Y'_{v'}=X'_{v'}.$ This procedure eliminates 
the particular pair of $u',\,v'\in 2^{m+1},$ as required. 

Then $X'=\bigcap_m\bigcup_{u\in 2^m}X_u$ belongs to $\pei.$ We 
prove that $X'$ is as required, that is, $F$ captures $\i$ on 
$X'.$ Assume that, on the contrary, there exists a pair of points 
$x,\,y\in X'$ such that $F(x)=F(y)$ but $x(\i)\not=y(\i).$ Let 
$x=x_a$ and $y=x_b$ (see the beginning of the proof of 
Theorem~\ref{fut}), $a,\,b\in 2^\om.$ Then 
$\i\not\in\xi[a,b]=\bigcap_m\xi[a\res m,b\res m]$ (see assertion 
$(\ast)$ in the proof of Theorem~\ref{fut}). Let $m$ be the least 
among those satisfying $\i\not\in\xi=\xi[a\res m,b\res m].$ Then 
$\xi\sq[\not\>\lsp\i],$ so that the case $(1)$ in $(\star)$ is 
impossible for $u=a\res m$ and $v=b\res m.$ Therefore 
$F\ima X_u\cap F\ima X_v=\emptyset,$ contradiction with the choice 
of $x$ and $y$ because $x\in X_u,$ $y\in X_v$.\qed\vspace{4mm}

We are already equipped enough to handle different items of 
Theorem~\ref{mm}.\vom

{\it Item \ref{mm2}\/}. Thus suppose that $F$ is reducible to 
$\xi$ on\/ $X$ and, on the contrary, $F$ does capture some 
$\i\in\bI\setminus\xi$ on $X.$ Then the co-ordinate function 
$C_\i(x)=x(\i)$ is itself reducible to $\xi$ on $X.$ Since 
$\i$ does not belong to $\xi,$ and on the other hand $C_\i$ is 
obviously reducible to $[\<\lsp\i],$ we conclude that $C_\i$ is also 
reducible to $[<\lsp\i]$ on $X$ by Lemma~\ref{inter} But this 
clearly contradicts requirement~\ref{perf1}.\vom\vom

{\it Item \ref{mm1+}\/}. Otherwise, by Lemma~\ref{l2} $C_\j$ is 
reducible to\/ $\xi=[\not\>\lsp\i]$ on some $X'\in\pei,$ 
$X'\sq X,$ contradiction with the already proved 
item~\ref{mm2}.\vom\vom

{\it Item \ref{mm2+}\/}. Arguing as in the proof of Lemma~\ref{l2}, 
we get a fusion system $\ang{X_u:u\in 2^{<\om}}$ such that 
$X_\La\sq X$ and $(\star)$ holds. We prove that the set 
$Y=\bigcap_n\bigcup_{u\in 2^m} X_u$ is as required. Suppose that on 
the contrary $x,\,y\in Y$ satisfy $F(x)=F(y);$ we have to prove that 
$x\res\xi=y\res\xi.$ By definition, $\ans{x}=\bigcap_n X_{a\res n}$ 
and $\ans{y}=\bigcap_n X_{b\res n}$ 
for certain (unique) $a,\,b\in 2^\om.$ It suffices to verify 
that $\xi\sq\xi[a\res m,b\res m]$ for all $m$.

Assume that on the contrary $\xi\not\sq\xi[u,v],$ where $u=a\res m$ 
and $v=b\res m$ for some $m.$ We assert that the case $(1)$ of 
$(\star)$ cannot happen. Indeed otherwise in particular $F$ is 
reducible to $\xi[u,v]$ on a set $X'=X_u\sq X.$ Take an arbitrary 
$\i\not\in\xi.$ By the assumption of item~\ref{mm2+} $F$ captures $\i$ 
on a set $X''\in\pei,$ $X''\sq X'.$ Thus the co-ordinate function 
$C_\i$ is reduced to $\xi[u,v]$ on $X''$ -- contradiction with the 
already proved item~\ref{mm2}. Thus we have case $(2)$ of $(\star),$ 
that is, $F\ima X_u\cap F\ima X_v=\emptyset.$ 
But this contradicts the assumption $F(x)=F(y)$.\vom\vom

{\it Item \ref{mm3}\/}. We fix $\Phi$ and put $\xi[u,v]=\xi_\Phi[u,v]$ 
as in the begining of the proof of Lemma~\ref{l2}. Assume that a set 
$X'\in\pei$ of type (b) of item~\ref{mm3} does not exist. 

Then by Lemma~\ref{l2} if 
$\i\in\bI\setminus\xi$ then every set $Y\in\pei,$ $Y\sq X$ 
contains a subset $Z\in \pei$ such that $F$ is reducible to 
$[\not\>\lsp\i]$ on $Z.$ Using lemmas \ref{suz} and \ref{pand} we 
obtain a fusion sequence $\ang{X_u:u\in 2^{<\om}}$ such that 
$X_\La\sq X$ and $F$ is reducible to $[\not\>\lsp\Phi(m)]$ on 
$X_u$ whenever $u\in 2^m$ and $\Phi(m)\not\in\xi.$ Then 
$X'=\bigcap_m\bigcup_{u\in 2^m}X_u\in\pei.$ We prove that $X'$ is 
a set of type (a), that is, $F$ is reducible to $\xi$ on $X'$.

Let us define $\xi_m\in\IS$ by induction on $m$ so that 
$\xi_0=\bI$ and 
\dm
\xi_{m+1}=\left\{
\bay{ccl}
\xi_m & - & \hbox{in the case }\, \Phi(m)\in\za;\\[3mm]

\xi_m\cap [\not\>\lsp\Phi(m)] & - & \hbox{in the case }\,
\Phi(m)\not\in\za
\eay
\right.
\dm
for all $m.$ Notice that then $\xi_m\sq\xi[u,v]$ whenever 
$u,\,v\in 2^m$ satisfy $\xi\sq\xi[u,v]$.\vspace{3mm}

\noi
{\bf Assertion} \ {\it For any\/ $m,$ $F$ is reducible to\/ 
$\xi_m$ on\/ $X_m=\bigcup_{u\in 2^m}X_u$.}\vspace{4mm}

\noi
{\bf Proof} \ of the assertion.  
We argue by induction on $m.$ Assume that, on the contrary, there 
exist $u,\,v\in 2^{m+1}$ and $x\in X_u,$ $y\in X_v$ such that 
$x\res\xi_{m+1}=y\res\xi_{m+1},$ but $F(x)\not= F(y).$ Let 
$u=u'\we d,$ $v=v'\we e,$ where $u',\,v'\in 2^m$ and 
$d,\,e\in\ans{0,1}.$ 

We have $\xi_{m+1}\sq\xi[u,v]$ by \ref{prct} and  
\ref{aprct}, therefore $\xi\sq\xi[u,v]$ because every set 
$\xi_n$ includes $\xi.$ This implies $\xi\sq\xi[u',v'].$ It 
follows (see above) that $\xi_m\sq\xi[u',v']$.

The nontrivial case is the case when $\i=\Phi(m)\not\in\xi$ since 
otherwise $\xi_{m+1}=\xi_m$ and we can use the inductive 
hypothesis. Then $\xi_{m+1}=\xi_m\cap [\not\>\i]$.

We assert that there exists $x'\in X_{u'}$ such that\its
\bit
\item $x'\res\xi_m=y\res\xi_m$\hspace{5mm} and\hspace{5mm}
$x'\rsd{\not\>\i}=x\rsd{\not\>\i}$.\its
\eit
Indeed, first $x\rsd{\not\>\i}\in X_{u'}\rsd{\not\>\i}.$ Second, 
since $\xi_m\sq\xi[u',v'],$ we obtain 
$X_{u'}\res\xi_m=X_{v'}\res\xi_m$ by \ref{prct}. Therefore 
$y\res\xi_m\in X_{u'}\res\xi_m.$ Finally, using the fact that 
$\xi_{m+1}=\xi_m\cap [\not\>\lsp\i]$ we conclude that 
$x\res(\xi_m\cap[\not\>\lsp\i])=y\res(\xi_m\cap[\not\>\lsp\i])$ 
by the choice of $x$ and $y.$ Property~\ref{perf1} of $X_{u'}$ 
then implies the existence of a point $x'\in X_{u'}$ satisfying 
$(\bullet)$.

To end the proof of the assertion, notice that $F(x')=F(y)$ by 
the inductive hypothesis while $F(x')=F(x)$ by the choice of the 
fusion sequence of sets $X_u$.\qed\vspace{4mm}

We continue the proof of item~\ref{mm3} of Theorem~\ref{mm}. 

It follows from the 
assertion that $F$ is reducible to every $\xi_m$ on $X'.$ This 
allows to conclude that $F$ is also reducible to $\xi$ on $X'.$ 
Indeed assume that on the contrary $x,\,y\in X'$ and 
$x\res\xi= y\res\xi$ but $F(x)\not=F(y).$ By the continuity of 
$F$ there exist $m\in\om$ and $u,\,v\in 2^m$ such that $x\in X_u,$ 
$y\in X_v,$ and $F\ima X_u\cap F\ima X_v=\emptyset.$ Notice that 
then ${X_u\res\xi\cap X_v\res\xi\not=\emptyset},$ therefore 
$\xi\sq\xi[u,v]$ by {\ref{aprct}}.  
This implies $\xi\sq\xi_m\sq\xi[u,v],$ as above. Therefore $F$ is 
reducible to $\xi[u,v]$ on $X',$ contradiction with the equality 
$F\ima X_u\cap F\ima X_v=\emptyset,$ because 
$X_u\res\xi[u,v]=X_v\res\xi[u,v]$ by condition \ref{prct}.\qed

\newpage
\subsection{Proof of the theorem: the countable case}
\label{re}

This section starts the proof of Theorem~\ref{m} in the case when the 
``length'' $\bI$ of iteration is countable. 
Thus let $\gM$ be a countable 
transitive model of $\ZFC,$ $\bI\in \gM$ be a countable in $\gM$ 
partially ordered set.\vspace{2mm}

\noi{\bf The forcing}\\[2mm] 
We consider $\peim=(\pei)^\gM$ as a forcing 
notion ($X\sq Y$ means that $X$ is a stronger condition). Notice 
that every set in $\dP$ is then a countable subset in the 
universe. However we can transform it to a perfect set by the 
closure operation: the topological closure $X^\#$ of a set 
$X\in \peim$ will then satisfy the definition of $\pei$ from the 
point of view of the universe.\vspace{2mm}

\noi{\bf The model}\\[2mm]
Let $G\sq\peim$ be a \dd\peim generic 
ultrafilter over $\gM.$ Then the intersection 
$\Pi=\bigcap_{X\in G}X^\#$ is a singleton (this easily follows 
from Proposition~\ref{clop}). Let $\x\in\cani$ be the unique 
element of $\Pi;$ thus $\x$ is a function from $\bI$ to reals. As 
usual in this case the generic extension $\gN=\gM[G]$ is equal 
to $\gM[\x]$.

We define $\a_\i=\x(\i)$ for all $\i\in\bI,$ so that 
$\x=\ang{\a_\i:\i\in\bI}$.

\bpro
\label{counta}
The model\/ $\gN=\gM[\x]$ satisfies the 
requirements of Theorem~\ref{m}.
\epro
The proof of this proposition is the content of this section. 
First we prove the cardinality 
preservation property and an important technical theorem which 
will allow to study reals in the extension using continuous 
functions in the initial model. 

\bte
\label{card}
$\aleph_1^\gM$ and any cardinal greater than\/ 
$2^{\aleph_0}$ in\/ $\gM$ remain cardinals in $\gN$.
\ete
\proof Let $\unf$ be a name of a function defined on $\om$ in the 
language of forcing. Using in $\gM$ lemmas \ref{suz} and 
\ref{pand} we obtain a fusion sequence $\ang{X_u:u\in 2^{<\om}}$ 
of sets $X_u\in\peim$ such that $X_\La\sq X_0$ (where $X_0\in\peim$ 
is a given condition) and for all $m,$ every $X_u,\,\,u\in 2^m,$ 
decides the value of $\unf(m).$ Then 
$X=\bigcap_m\bigcup_{u\in 2^m}X_u\in\peim,$ $X\sq X_0,$ and $X$ 
forces that the range of $\unf$ is a subset of a countable in 
$\gM$ set.\qed\vspace{4mm}

\noi{\bf Continuous functions}\\[2mm]
We put $\cntm\xi=(\cnt\xi)^\gM$ and $\cntim=(\cnti)^\gM.$ 
It is a principal property of several forcing notions 
(including Sacks forcing) that reals in the generic extension can 
be obtained by application of continuous functions (having a code) 
in the ground model, to the generic object. As we shall prove this 
is also a property of the Sacks iteration considered here. 

Obviously every $F\in\cntm\xi$ is a countable subset 
of $\can\xi\ti \om^\om$ in the universe, but since the domain of 
$F$ in $\gM$ is the compact set $\can\xi,$ $F^\#$ is a continuous 
function mapping $\can\xi$ into the reals in the universe. 

\bte
\label{repr}
Let\/ $\xi\in\IS$ and $c\in\gM[\x\res\xi]\cap\om^\om.$ There 
exists a function\/ $H\in\cntm\xi$ such that\/ 
$c=H^\#(\x\res\xi)$.~\footnote
{\rmt\ Obviously this equality is absolute for any model containing 
all of $c,\,\x,\,H$.} 
\ete
\proof Let $\unc$ be a name for $c$ containing an explicit 
absolute construction of $c$ from $\x\res\xi$ and some parameter 
$p\in\gM.$  

{\it We argue in $\gM.$} 

It follows from lemmas \ref{pro'} and 
\ref{apro} that forcing of properties of $\unc$ is reduced to 
$\xi$ in the sense that if $X$ forces $\unc(m)=k$ then 
$X^-=X\res\xi\ares\bI$ forces $\unc(m)=k,$ too. 
Therefore, given $X_0\in\pei,$ we can define a fusion sequence 
$\ang{X_u:u\in 2^{<\om}}$ of sets $X_u\in\pei$ such that 
$X_\La\sq X_0$ and, for all $m\in\om$ and $u\in 2^m,$ $X_u$ 
decides the value of $\unc(m),$ say, forces that 
$\unc(m)=\al(u),$ where $\al\in\gM$ maps $2^{<\om}$ into $\om,$ 
so that if $X_u\res\xi=X_v\res\xi$ then $\al(u)=\al(v).$ The 
function $F'$ defined on 
$X=\bigcap_m\bigcup_{u\in 2^m}X_u\in\pei$ by the property: 
$F'(x)(m)=k$ iff $\al(u)=k,$ -- for all $m$ and all $u\in 2^m$ 
such that $x\in X_u,$ is continuous and reducible to $\xi$ on 
$X.$ Therefore $F'$ can be expanded to a function $F\in\cnti$ 
reducible to $\xi$ on $\cani.$ In this case (see 
Remark~\ref{redr}) there exists a function $H\in\cnt\xi$ such 
that $F(x)=H(x\res\xi)$ for all $x\in\cani.$ Then $X$ forces that 
$\unc=F^\#(\x)=H^\#(\x\res\xi)$.\qed\vspace{3mm}

This theorem practically reduces properties of reals in 
\dd\peim generic extensions to properties of continuous functions 
in the ground model.

To demonstrate how Theorem~\ref{repr} works we prove statements 
\ref{m1+} through \ref{m3} of Theorem~\ref{m} for the model 
$\gN=\gM[G]$ of consideration.\vspace{3mm}

\noi{\bf Proof of statement \ref{m1+} of Theorem~\ref{m}}\\[2mm]
Thus let $\i,\,\j\in\bI,\;\;\i<\j;$ we have to prove that 
$\a_\i\in\gM[\a_\j].$ It follows from Theorem~\ref{mm} 
(item~\ref{mm1+}) that there exists a condition $X\in G$ such that 
the following is true in $\gM:$ the co-ordinate function $C_\j$ 
captures $\i$ on $X.$ In other words, in $\gM$ there exists a 
continuous function $H:\cD\,\lra\,\cD$ such that 
$x(\i)=H(x(\j))$ for all $x\in X.$ It follows that 
$x(\i)=H^\#(x(\j))$ for all $x\in X^\#$ is also true in $\gN.$ 
($H^\#$ is the topological closure of $H$ as a subset of $\cD^2$.) 
Therefore $\a_\i=H^\#(\a_\j)\in \gN=\gM[G]$.\vspace{3mm}

\noi{\bf Proof of statement \ref{m2} of Theorem~\ref{m}}\\[2mm]
Let $\xi\in \gM$ be an 
initial segment of $\bI,$ and $\i\in\bI\setminus\xi;$ we have 
to prove that $\a_\i\not\in\gM[\x\res\xi].$ Assume that on the 
contrary $\a_\i\in\gM[\x\res\xi].$ Theorem~\ref{repr} implies 
$\a_\i=H^\#(\x\res\xi)$ for a function $H\in\cntm\xi.$ Let this 
be forced by some $X\in\peim$ such that $\x\in X^\#$.

{\it We argue in $\gM$.} 

The set 
${Y=\ans{x\in X:x(\i)=H(x\res\xi)}}$ is a closed subset of $X$ 
because $H$ is continuous. Therefore by Corollary~\ref{bor} 
there exists $X'\in\pei$ such that either $X'\sq Y$ or 
$X'\sq X\setminus Y.$ The former possibility means that the 
function $F(x)=x(\i)$ is reducible to $\xi$ on $X'.$ Since $F$ 
is obviously reducible to $[\<\lsp\i],$ it is also reducible to 
$\xi\cap[\<\lsp\i]$ by Lemma~\ref{inter}, therefore to 
$[<\lsp\i],$ because $\i\not\in\xi.$ But this evidently 
contradicts property~\ref{perf1} of $X'$.

This contradiction shows that in fact $X'\sq X\setminus Y.$ But 
then $X'$ forces that $\x(\i)=\a_\i$ is {\it not\/} equal to 
$H^\#(\x\res\xi),$ contradiction with the choice of $X$ and 
$X'$.\qed\vspace{3mm}

\noi{\bf Proof of statement \ref{m2+} of Theorem~\ref{m}}\\[2mm]
We obtain the result from item~\ref{mm2+} of Theorem~\ref{mm} 
using essentially the same type of reasoning as above.\vspace{3mm}

\noi{\bf Proof of statement \ref{m3} of Theorem~\ref{m}}\\[2mm]
Let $\xi\in \gM$ be an 
initial segment of $\bI,$ and $c\in\gN\cap\cD.$ We have to prove 
that either $c\in\gM[\x\res\xi]$ or there exists $\i\not\in\xi$ 
such that $\a_\i\in\gM[c].$ Theorem~\ref{repr} tells that 
$c=F^\#(\x)$ for some $F\in\cntim.$ Let this be forced by some 
$X\in\peim.$ We also assume that on the contrary $c$ 
{\it does not\/} satisfy requirement~\ref{m3} of 
Theorem~\ref{m}, and this is forced by the same $X$.

{\it We argue in $\gM$.} 

It follows from Theorem~\ref{mm} (item~\ref{mm3}) that 
there exists $X'\in\pei,$ $X'\sq X,$ such that either $F$ is 
reducible to $\xi$ on $X'$ or $F$ captures some $\i\not\in\xi$ 
on $X'.$ 

Consider the ``either'' case. Then (see Remark~\ref{redr}) there 
exists a function $H\in\cnt\xi$ such that $F(x)=H(x\res\xi)$ for 
all $x\in X'.$ In this case $X'$ forces that $c\in\gM[\x\res\xi],$ 
contradiction with the choice of $X$ and $X'$.

Consider the ``or'' case. Then there exists a continuous function 
$E:\cD\;\lra\;\cD$ such that $x(\i)=E(F(x))$ for all $x\in X'.$ 
Then $X'$ forces $\a_\i=\x(\i)=E^\#(F^\#(\x))\in\gM[F^\#(\x)],$ 
again a contradiction with the choice of $X$ and $X'$.

{\it We argue in $\gN$.} 

It remains to note that the possibilities 
of statement~\ref{m3} of Theorem~\ref{m} are incompatible by the 
already proved statement~\ref{m2}.\qed

\subsubsection*{The Sacksness}

In this subsection we prove the principal item of Theorem~\ref{m} 
--- statement~\ref{m1} which shows that in fact $\gN$ is an 
iterated Sacks extension of $\gM$ with $\bI$ as the ``length'' of 
the iteration. 

\bte
\label{sack}
Every\/ $\a_\i$ is Sacks generic over 
$\gM[\x\rsd{<\i}]=\gM[\ang{\a_\j:\j<\i}]$. 
\ete
\proof Assume that $S\in\gM[\x\rsd{<\i}]$ is, in 
$\gM[\x\rsd{<\i}],$ a dense subset in the collection of all 
perfect subsets of $\cD=2^\om;$ we have to prove that 
$\a_\i\in C^\#$ for some $C\in S.$ Suppose that on the contrary 
some $X\in \peim$ such that $\x\in X^\#$ forces the opposite. 

{\it We argue in $\gM.$\/} 

The set 
$D(y)=D_{Xy}(\i)=\ans{x(\i):x\in X\cj x\rsd{<\i}=y}$ is perfect 
for all $y\in Y=X\rsd{<\i}$ by property~\ref{perf1} of $X$. 

{\it We argue in $\gM[\x\rsd{<\i}].$\/} 

Notice that 
$\y=\x\rsd{<\i}$ belongs to $Y^\#.$ Therefore $D^\#(\y)$ is a 
perfect set (certain absoluteness is applied). Thus there exists 
a set $C\in S$ such that $C\sq D^\#(\y)$.

Now we are in need of a coding of closed subsets of $\cD=2^\om.$ 
Let $\ans{B_n:n\in\om}$ be a recursive enumeration of 
all basic clopen sets in $\cD.$ We put 
$\clo{c}=\cD\setminus\bigcup_{c(n)=0}B_n$ for all $c\in 2^\om.$ 
Thus every closed $C\sq\cD$ is equal to $\clo c$ for some 
$c\in 2^\om$. 

{\it We argue in $\gM[\x\rsd{<\i}]$}. 

We define $c\in 2^\om$ so 
that $c(n)=0$ iff $B_n$ is a basic clopen set disjoint from $C;$ 
then $C=\clo{c}.$ By Theorem~\ref{repr}, $c=H^\#(\y)$ for 
some $H\in\cntm{<\i}.$ Since $\a_\i=\x(\i)\not\in C,$ we can 
assume that $X$ forces that $\clo{H^\#(\x\rsd{<\i})}$ is a perfect 
subset of $D^\#(\x\rsd{<\i})$ and $\x(\i)$ does not belong to 
$\clo{H^\#(\x\rsd{<\i})}$. 

{\it We argue in $\gM.$\/} 

The continuation of the proof is 
based on the following fact: the set
\dm 
U=\ans{x\in X: x(\i)\in\clo{H(x\rsd{<\i})}}
\dm
contains a subset 
$U'\in\pei.$ Such a condition $U'$ would force that $\x(\i)=\a_\i$ 
belongs to $\clo{H^\#(\x\rsd{<\i})}$ by absoluteness, 
contradiction with the statement forced by $X$. 

Thus we concentrate on the mentioned key fact. Notice that 
since we deal with compact spaces, the set 
\dm
Y_1=\ans{y\in Y:\clo{H(y)}\,\hbox{ is a perfect subset of }\,D(y)}
\dm
is $\bbox{G}_\da$ in $Y.$ Corollary~\ref{bor} says that there 
exists a set $Y_2\in\pe{<\i},$ $Y_2\sq Y,$ such that either 
$Y_2\sq Y_1$ or $Y_2\cap Y_1=\emptyset.$ Then 
$X_2=X\cap (Y_2\ares\bI)\in\pei$ by Lemma~\ref{apro}. 

Suppose that $Y_2\cap Y_1=\emptyset.$ Then by absoluteness $X_2$ 
forces that $\clo{H^\#(\x\rsd{<i})}$ is not a perfect subset of 
$D^\#(\x\rsd{<i}),$ contradiction with the statement forced by $X$. 

Thus in fact $Y_2\sq Y_1.$ We have to restrict $Y_2$ a little bit 
more. Notice that for any clopen $G\sq\cD$ the set 
$Y(G)=\ans{y\in Y_2:\clo{H(y)}\cap G\not=\emptyset}$ is closed. 
Then by Corollary~\ref{prebor} there exists $Y'\sq Y_2,$ 
$Y'\in\pe{<\i}$ such that $Y'\cap Y(G)$ is clopen in $Y'$ 
for every $G$. 

We demonstrate that 
${Z=\ans{z\in X\rsd{\<\i}:z\rsd{<\i}\in Y'\cj z(\i)\in 
\clo{H(z\rsd{<\i})}}}$ 
belongs to $\pe{\<\i}.$ To check all the requirements of 
lemmas \ref{gath} and \ref{gath'} notice that $Z$ is closed and 
$Z\rsd{<\i}=Y'\in\pe{<\i}.$ Moreover if $y\in Y'$ then 
$D_{Zy}(\i)=\clo{H(y)}$ is perfect since $Y'\sq Y_1,$ so that 
$Z$ satisfies \ref{perf1}. Finally $Z$ also satisfies \ref{oz} by 
the choice of $Y'.$ Therefore in fact $Z\in\pe{\<\i}$.

It remains to set $U'=X\cap (Z\ares\bI)$.\qed\vspace{4mm}

This ends the proof of Proposition~\ref{counta} --- the countable 
case in Theorem~\ref{m}.\qed

\newpage
\subsection{Uncountable case}
\label{uncount}

We carry out the general case of Theorem~\ref{m} in very brief 
manner because the principal points can be reduced to the already 
considered countable case. 

Thus let $\gM$ be a countable transitive model of $\ZFC,$ 
$\bI\in\gM$ a po set. 

Let $\cs$ be the collection of all sets $\xi\in\gM,$ $\xi\sq\bI,$ 
such that $\card\xi\<\aleph_0$ in $\gM.$ Notice that sets 
$\xi\in\cs$ are, generally speaking, not initial segments of 
$\bI$ or of each other. 

For any $\xi\in\cs,$ let $\peM\xi=(\pe\xi)^\gM.$ The set  
$\peim=\bigcup_{\xi\in\cs}\peM\xi$ is the forcing notion. To 
define the order, we first put $\supp X=\xi$ whenever 
$X\in\peM\xi.$ Now we set $X\<Y$ ($X$ is stronger than $Y$) iff 
$\xi=\supp Y\sq \supp X$ and $X\res\xi\sq Y$. 

Let $G\sq\peim$ be a generic set over $\gM.$ Then there exists 
unique indexed set $\x=\ang{\a_\i:\i\in\bI},$ all $\a_\i$ belong 
to $\cD,$ such that $\x\res\xi\in X^\#$ whenever $X\in G$ and 
$\supp X=\xi.$ Moreover $\gM[G]=\gM[\x]=\gM[\ang{\a_\i:\i\in\bI}]$.

\bpro
\label{m-m} 
The model\/ $\gN=\gM[G]=\gM[\x]$ satisfies Theorem~\ref{m}. 
\epro
\proof{}is based on the two principal statements. 

\bte
\label{gen1}
$\aleph_1^\gM$ remains a cardinal in $\gN$.~\footnote
{\rmt\ The behaviour of other cardinals depends on the cardinal 
structure in $\gM,$ the cardinality of $\bI,$ and the cardinality 
of chains in $\bI.$ It is not our intension here to investigate 
this matter.}
\ete

\bte
\label{gen2}
Assume that\/ $\bJ\in\gM$ is an initial segment of\/ $\bI$ and\/ 
$c\in\gM[\x\res\bJ]\cap\om^\om.$ There exist\/ 
$\xi\in\cs,$ $\xi\sq\bJ,$ and a function\/ $H\in\cntm\xi$ such 
that\/ $c=H^\#(\x\res\xi)$. 
\ete

Theorem~\ref{gen2} allows to repeat the reasoning in 
Section \ref{re} and prove statements \ref{m1} through \ref{m3} 
of Theorem~\ref{m} using properties of forcing conditions and 
continuous functions proved mainly in Section~\ref{cont}. 
Thus theorems \ref{gen1} 
and \ref{gen2} suffice for Proposition~\ref{m-m} and 
Theorem~\ref{m}.\qed\qed\vspace{3mm}

\noi
{\bf Proof}\hspace{1mm} of Theorem~\ref{gen1}\\[2mm]
Let $\unf$ be a name of a function defined on $\om$ in the 
language of forcing. We fix $X_0\in\peim.$ The aim is to 
obtain a condition $X\in\peim,$ $X\<X_0,$ and a countable in 
$\gM$ set $R$ such that $X$ forces that the range of $\unf$ is 
included in $R.$ Let $\xi_0=\supp{X_0}$.\vom

{\it We argue un $\gM$.}
 
To utilize the proof of Theorem~\ref{card} we reduce the forcing 
of statements related to $\unf$ to a certain $\za\in\cs.$ 

Let $\xi\in\cs.$ We say that a set $\cX\sq\pe\za$ is {\it 
adequate\/} if\its 
\ben
\def\theenumi{\alph{enumi})}
\def\labelenumi{\theenumi}
\itla{a}
for any initial segment 
$\eta=
\ans{\j\in\xi: \j\not\>\i_1\cj...\cj \j\not\>\i_n},$ 
where $\i_1,...,\i_n\in\xi,$ 
and any pair $X,\,Y\in\cX,$ if $Y\res\eta\sq X\res\eta$ then 
$Z=X\cap (Y\res\xi\ares\eta)\in\cX,$ \hfill and\hfill\its

\itla{b} 
for all $X\in\cX,$ $\i\in\xi,$ $e\in\ans{0,1},$ the 
set $X_e=\spl(X,\i,e)$ belongs to $\cX.$~\footnote
{\rmt\ Notice that $Z\in\pe\xi$ and $X_e\in\pe\xi$ by 
lemmas \ref{apro} and \ref{spl}.}\its 
\een
It is obvious that every countable $\cX'\sq\pe\xi$ can be extended 
to a countable adequate $\cX\sq\pe\xi$.

It can be easily verified that if $\xi\sq\za\in\cs$ and 
$X\in\pe\xi$ then $X\ares\za\in\pe\xi$ (although in general 
Lemma~\ref{apro} is not true in the case when $\xi$ is not an 
initial segment of $\za$). Therefore for all $\xi\in\cs,$ 
$n\in\om,$ and countable $\cX'\sq\pe\xi$ there exists $\za\in\cs$ 
and an adequate countable $\cX\sq\pe\za$ such that $\xi\sq\za$ 
and\its
\ben
\def\theenumi{\roman{enumi})}
\def\labelenumi{\theenumi}
\item 
$X'\ares\za\in \cX$ whenever $X'\in\cX'$;\hfill and\hfill\its

\item
for any $X'\in\cX'$ there exists $X\in\cX$ such that $X\<X'$ 
and $X$ decides the value of $\unf(n)$.\its
\een
This allows to start from $X_0\in\pe{\xi_0}$ and define by 
induction a sequence $\xi_0\sq\xi_1\sq\xi_2\sq...$ of 
$\xi_n\in\cs$ and a sequence of countable adequate 
$\cX_n\sq\pe{\xi_n}$ such that\its
\ben
\def\theenumi{\arabic{enumi})}
\def\labelenumi{\theenumi}
\item 
$X\ares\xi_{n+1}\in \cX_{n+1}$ whenever $X\in\cX_n$;\hfill 
and\hfill\its

\item 
for any $X'\in\cX_n$ there exists $X\in\cX_{n+1}$ such that 
$X\<X'$ and $X$ decides the value of $\unf(n)$.\its
\een
We set $\za=\bigcup_{n\in\om}\xi_n$ and 
$\cX=\bigcup_{n\in\om}\ans{X\ares\za:X\in\cX_n}.$ Then 
$\cX\sq\pe\za$ is a countable adequate family which satisfies the 
property that 
\dm
\forall\,X'\in\cX\;\forall\,n\;\exists\,X\in\cX\;
[\,X\sq X'\cj X\,\hbox{ decides the value of }\,\unf(n)\,]\,.
\dm

We notice now that the transformations of sets used in the proofs 
of lemmas \ref{suz} and \ref{pand} are of types \ref{a} and 
\ref{b}. Therefore arguing like in the proof of Theorem~\ref{card} 
we can obtain a fusion sequence $\ang{X_u:u\in 2^{<\om}}$ 
of sets $X_u\in\cX$ such that $X_\La\sq X_0$ and for all $m,$ 
every $X_u,\,\,u\in 2^m,$ decides the value of $\unf(m).$ Then 
$X=\bigcap_m\bigcup_{u\in 2^m}X_u\in\pe\za,$ $X\sq X_0,$ and $X$ 
forces that the range of $\unf$ is a subset of a countable in 
$\gM$ set.\qed\vspace{3mm}

\noi
{\bf Proof}\hspace{1mm} of Theorem~\ref{gen1}\\[2mm]
Let $\unc$ be a 
name for $c$ containing an explicit absolute construction of $c$ 
from $\x\res\bJ$ and some $p\in\gM$.  

{\it We argue in $\gM.$} 

Given $X_0\in\pe{\xi_0},$ we argue as in 
the proof of Theorem~\ref{gen1} and get $\za\in\cs,$ 
$\xi_0\sq\za\sq\bJ,$ and a countable adequate $\cX\sq\pe\za$ 
such that $X_0\ares\za\in\cX$ and 
\dm
\forall\,X'\in\cX\;\forall\,n\;\exists\,X\in\cX\;
[\,X\sq X'\cj X\,\hbox{ decides the value of }\,\unc(n)\,]\,.
\dm
It remains to carry out the construction in the proof of 
Theorem~\ref{repr} within $\cX$.\qed\vspace{3mm}

This ends the proof of Theorem~\ref{m} in general 
case.\qed

\newpage 
\newcommand{\xxi}{{\hspace{1pt}\underline{\x\res\xi}\hspace{1pt}}}
\newcommand{\xeta}{{\hspace{1pt}\underline{\x\res\eta}\hspace{1pt}}}
\subsection{The non--Glimm--Effros equivalence 
relation}
\label{seqge}

This section is devoted entirely to the proof of Theorem~\ref{ge}. 
Thus let $\C$ be the equivalence relation defined on reals by 
$x\C y$ iff $\rL[x]=\rL[y]$. 

We have to find a model for Theorem~\ref{ge}. 
Let $\gM$ be a countable transitive model of $\ZFC$ plus the 
axiom of constructibility -- the initial model. 

We define the ``length'' of the iteration $\bI=\om_1^\gM\ti\dZ$ 
($\om_1^\gM$ copies of $\dZ,$ the integers). Thus from the point 
of view of $\gM,$ $\bI$ is the set of all pairs $\ang{\al,z},$ 
$\al<\om_1$ and $z\in \dZ,$ linearly ordered lexicographically, 
but of course not wellordered. 

Let $\cs$ be the collection of all initial segments $\xi\in\gM,$ 
$\xi\sq\bI,$ such that $\card\xi\<\aleph_0$ in $\gM.$ (This 
formally differs from the definition in Section~\ref{uncount}, 
but not essentially, since each \dd\gM countable subset of $\bI$ 
can be covered by a countable initial segment.)

We define $\peM\xi=(\pe\xi)^\gM$ (for all $\xi\in\cs$), 
$\peim=\bigcup_{\xi\in\cs}\peM\xi,$ the ``support'' $\supp X,$ 
and the order on $\peim$ as in Section~\ref{uncount}. 

Let us fix a generic over $\gM$ set $G\sq\peim.$ We define 
$\x=\ang{\a_\i:\i\in\bI}$ (all $\a_\i$ being elements of $\cD$). 
We prove that the model 
$\gN=\gM[G]=\gM[\x]=\gM[\ang{\a_\i:\i\in\bI}]$ satisfies 
Theorem~\ref{ge}. 

\bte
\label{tge}
It is true in\/ $\gN$ that\/ $\C\;:$\its
\bit
\item[--] neither admits a R-OD separating family$;$\its

\item[--] nor admits an uncountable R-OD pairwise\/ 
\dd\C inequivalent set$.$
\eit
\ete
\proof Let us first investigate the structure of the degrees of 
constructibility (i.e. \dd\C degrees) in $\gN$ -- or, that is the same 
since $\gM$ models $\rV=\rL,$ degrees of \dd\gM constructibility. 

The set 
$\bI$ has a nice property: its initial segments admit a clear description. 
Indeed, each $\xi\in\cs$ is equal to one of the following:
\dm
\bay{rccccl}
\xi_{\al z} & = & [\<\lsp\i] & = &\ans{\j\in\bI:\j\<\i}\,, & 
\hbox{where }\,\i=\ang{\al,z}\in\bI\,\hbox{ and }\,
\al<\om_1^\gM,\;\,z\in\dZ\,;\\[2mm]

\xi_{\al} &&& = & \al\ti\dZ & \hbox{for some }\,\al<\om_1^\gM\,;
\eay
\dm
obviously all of them belong to $\gM.$ In particular, $\cs\in\gM$ is 
linearly ordered in $\gM$ by inclusion. One can see that $\cs$ is order 
isomorphic to $\om_1^\gM\ti(1+\dZ)$.

\ble
\label{ndeg}
Suppose that\/ $c$ is a real in\/ $\gN.$ There exists unique\/ 
$\xi\in\cs$ such that\/ $\gM[c]=\gM[\x\res\xi]$.
\ele
\proof{}of the lemma. (We recall that $\x\res\xi=\ang{\a_\i:\i\in\xi}$.) 
By Theorem~\ref{gen2}, $c=H^\#(\x\res\za)$ for appropriate $\za\in\cs$ 
and $H\in\cnt\za$ in $\gM;$ in particular $c\in\gM[\x\res\za]$. 

We set $\xi=\ans{\i\in\za:\a_\i\in\gM[c]};$ then $\xi\in\cs$ by 
Theorem~\ref{m}, item~\ref{m1+}. (We mean: by Proposition~\ref{m-m}  
which assures that $\gN$ satisfies 
Theorem~\ref{m}.) Furthermore items \ref{m2}, \ref{m2+}, 
\ref{m3} of the same theorem prove that both $c\in\gM[\x\res\xi]$ 
and $\x\res\xi\in\gM[c]$.\qed

\bcor
\label{cmany}
Suppose that\/ $c$ is a real in\/ $\gN.$ There exists only 
countably\/ {\rmt(in $\gN$)} many\/ \dd\gM degrees below $c$.\qed
\ecor

\subsubsection*{Proof of the ``nor'' part of Theorem~\ref{tge}}
 
Let, in $\gN,$ $S$ be a \dd\C pairwise inequivalent subset of $\cD$ 
defined (in $\gN$) by a formula containing ordinals and a real 
$p\in\gN$ as parameters. By Lemma~\ref{ndeg} there exists an 
initial segment, say $\eta=\xi_\al=\al\ti\dZ,$ $\al<\om_1^\gM,$ 
such that $p\in\gM[\x\res\eta].$ Then $S$ is definable in $\gN$ 
using a formula containing $\x\res\eta=\ang{\a_\i:\i\in\eta}$ and 
ordinals as parameters.\its 
\bit
\item[$(1)$] {\it We assert that\/ $S\sq\gM[\x\res\eta]$}.\its
\eit
This assertion implies that $S$ is countable by 
Corollary~\ref{cmany}, therefore suffices for the ``nor'' part of 
the theorem.

To prove the assertion, let us fix a real $r\in S$ in $\gN.$ Then 
we have $\rL[r]=\rL[\x\res\xi]$ for certain (unique) $\xi\in\cs$ by 
Lemma~\ref{ndeg}. 
Let $\vpi(\x\res\eta,\x\res\xi,k,l)$ be the formula:\its
\bit
\item there exists a real $r'\in S$ such that 
$\rL[r']=\rL[\x\res\xi]$ and $r'(k)=l$.\its
\eit
($\x\res\eta$ enters the formula via a definition of $S.$ We recall 
that $\rV=\rL$ is assumed in $\gM,$ so that $\rL[...]$ in $\gN$ 
corresponds to $\gM[...]$ in the universe.) Then, 
in $\gN,$ we have $r(k)=l$ iff $\vpi(\x\res\eta,\x\res\xi,k,l)$ for 
all $k,\,l$.

Let $\xxi$ and $\xeta$ be the names for $\x\res\xi$ and 
$\x\res\eta$.\its
\bit
\item[$(2)$] {\it We assert that, for all\/ $k,\,l\in\om$ and\/ 
$X\in\pei,$ if $\eta\sq\supp X$ and\/ $X$ forces\/ 
$\vpi(\xeta,\xxi,k,l)$ then\/ $X\res\eta$ already forces 
$\vpi(\xeta,\xxi,k,l)$.}\its
\eit
One can easily see that $(2)$ implies $r\in\gM[\x\res\eta],$ that 
is, implies $(1);$ therefore we concentrate on the assertion $(2)$. 

Assume that $(2)$ is not true. Thus results in a pair of conditions 
$X,\,Y\in \peM\za,$ where $\za\in\cs,$ $\eta\sq\za,$ such that 
$X\res\eta=Y\res\eta,$ $X$ forces $\vpi(\xeta,\xxi,k,l),$ but $Y$ 
forces the negation of $\vpi(\xeta,\xxi,k,l),$ for some $k,\,l$.

In $\gM,$ both $X$ and $Y$ are members of $\pe\za.$ Let $F\in\gM$ 
be a homeomorphism $X$ onto $Y$ satisfying requirements \ref{h1}, 
\ref{h2}, \ref{h3} of Theorem~\ref{thom}, in particular, 
$x\res\eta=F(x)\res\eta$ for all\/ $x\in X,$ because 
$X\res\eta=Y\res\eta$.

(Let us forget temporarily that a generic set $G\sq\peim$ was fixed 
above.) The homeomorphism $F$ induces the total automorphism of 
the part of $\peim$ stronger than $X$ onto 
the part of $\peim$ stronger than $Y,$ which results in a pair of 
\dd\peim generic over $\gM$ sets $G,\,G'\sq\peim$ and corresponding 
$\x,\,\x'\in\cani$ such that $X\in G,$ $Y\in G',$ $\gM[G]=\gM[G'],$ 
$\x\res\eta=\x'\res\eta,$ and finally 
$\gM[\x\res\xi]=\gM[\x'\res\xi]$ for all $\xi\in\cs.$ Thus we have 
got one and the same generic extension $\gN=\gM[G]=\gM[G']$ 
using two different generic sets.

Notice that the statement $\vpi(\x\res\eta,\x\res\xi,k,l)$ is true 
while $\vpi(\x\res\eta,\x'\res\xi,k,l)$ is false in $\gN$ by the 
choice of $X,\,Y.$ We cannot assert that $\x\res\xi=\x'\res\xi$ 
(unless $\xi\sq\eta,$ of course), but the formula $\vpi$ was 
defined so that it is \dd\C invariant on the argument $\x\res\xi:$ 
in other words,
\dm
\vpi(\x\res\eta,\x\res\xi,k,l)\;\;\llra\;\;
\vpi(\x\res\eta,\x'\res\xi,k,l)
\dm
provided $\gM[\x\res\xi]=\gM[\x'\res\xi].$ Since this assumption 
was obtained above, we conclude that 
$\vpi(\x\res\eta,\x'\res\xi,k,l)$ must be true in $\gN,$ 
contradiction. 

This ends the proof of the ``nor'' part of Theorem~\ref{tge}. 
 
\subsubsection*{Proof of the ``neither'' part of 
Theorem~\ref{tge}}

Suppose that, on the contrary, $\C$ admits in $\gN$ a R-OD 
separating family $\ang{X_\al:\al<\ga},$ $\ga$ an ordinal. As 
above, then the family is definable by a formula containing 
ordinals and some $\x\res\eta,$ $\eta\in\cs,$ as parameters. We 
define, in $\gN$,
\dm
U(r)=\ans{\al<\ga:r\in X_\al}
\dm
for each real $r\in\gN.$ Thus, $x\C y$ iff $U(x)=U(y)$ for each pair 
of reals $x,\,y$ in $\gN$.\its 
\bit
\item[$(3)$] {\it We assert that\/ $U(r)\in\gM[\x\res\eta]$ for 
all reals $r$ in $\gN$}.\its
\eit
Generally speaking, one would expect that $U(r)$ needs $r,$ or 
at least the \dd\C degree of $r$ as a parameter of definition. 
However, the \dd\C degrees in $\gN$ form a quite regular structure 
by Lemma~\ref{ndeg}, so that 
each degree is ``almost'' ordinal definable (but actually 
{\it not\/} OD), which makes it possible to prove $(3)$. 

As before, in the proof of the ``nor'' part, we reduce $(3)$ to a 
forcing assertion. Let us fix a real $r\in\gN;$ then by 
Lemma~\ref{ndeg} there exists $\xi\in\cs$ such that 
$\rL[r]=\rL[\x\res\xi].$ Let $\vpi(\x\res\eta,\x\res\xi,\al)$ 
be the formula:\its
\bit
\item there exists a real $r'$ such that 
$\rL[r']=\rL[\x\res\xi]$ and $r'\in X_\al$.\its
\eit
($\x\res\eta$ enters the formula via the enumeration of sets 
$X_\al$.) Then, in $\gN,$ we have $\al\in U(r)$ iff 
$\vpi(\x\res\eta,\x\res\xi,\al)$ for all $\al$.\its
\bit
\item[$(4)$] {\it We assert that, for all\/ $\al<\ga$ and\/ 
$X\in\pei,$ if $\eta\sq\supp X$ and\/ $X$ forces\/ 
$\vpi(\xeta,\xxi,\al)$ then\/ $X\res\eta$ already forces 
$\vpi(\xeta,\xxi,\al)$.}\its
\eit
As in the proof of the ``nor'' part above, $(4)$ implies 
$U(r)\in\gM[\x\res\eta],$ that 
is, implies $(3);$ therefore it suffices to prove $(4).$ We omit the 
reasoning because it is a copy of the proof of $(2)$ above: the 
principal point is that the formula $\vpi$ is again \dd\C invariant 
on the argument $\x\res\xi$. 

Thus we obtain $(4)$ and $(3)$.

We continue the proof of the ``either'' part. The key moment is as 
follows. It follows from assertion $(3)$ that in $\gN$ each \dd\C 
degree is definable by a formula using only ordinals and $\x\res\eta$ 
as parameters. In particular, $\x\res\eta$ plus ordinals as parameters 
is enough to distinguish all \dd\C degrees from each other. 

This will help us to engineer a contradiction. The special mechanism  
of getting a contradiction is based on the existence of order 
automorphisms (shiftings) in each \dd\dZ group in 
$\bI=\om_1^\gM\ti\dZ$.

There exists an ordinal $\al<\om_1^\gM$ such that 
$\eta\sq\xi_\al=\al\ti\dZ.$ We set $\i=\ang{\al,0},$ 
$\i'=\ang{\al,1}$ -- two neighbouring elements in the least \dd\dZ 
group not participating in $\xi_\al.$ We set $\xi=[\<\lsp\i],$ 
$\xi'=[\<\lsp\i'].$ Since $\i'\in\xi'\setminus\xi,$ we have 
$\gM[\x\res\xi]\not=\gM[\x\res\xi'].$ (Item \ref{m2} of 
Theorem~\ref{m} via Proposition~\ref{m-m}.) Take a pair of reals 
$r,\,r'\in\gN$ such that $\gM[r]=\gM[\x\res\xi]$ and 
$\gM[r']=\gM[\x\res\xi'];$ then $\rL[r]\not=\rL[r']$ in $\gN,$ hence  
$U(r)\not=U(r').$ Since both $U(r)$ and $U(r')$ belong to 
$\rL[\x\res\eta]$ in $\gN$ by $(3),$ we conclude that there exists a 
formula $\psi(\x\res\eta,x)$ containing only ordinals and 
$\x\res\eta$ as parameters, and such that the following is true 
in $\gN$ for every real $x$:
\dm
\rL[x]=\rL[\x\res\xi]\;\lra\;\psi(\x\res\eta,x))
\hspace{8mm}\hbox{and}\hspace{8mm}
\rL[x]=\rL[\x\res\xi']\;\lra\;\neg\;\psi(\x\res\eta,x))\,.
\eqno{(\ast)}
\dm
Therefore, a condition $X\in G$ forces $(\ast).$ One can w.l.o.g. 
assume that $\eta\sq \supp X.$\its
\bit
\item[$(5)$] {\it We assert that the weaker condition 
$Y=X\res\eta$ forces $(\ast)$.}\its
\eit
This is an assrtion of the same type as $(2)$ and $(4)$ above; its 
proof does not differ from the proof of $(2)$. 

Notice that $Y\in G$.

Let us consider the order automorphism $h:\bI\,\hbox{ onto }\,\bI$ 
defined as follows: $h(\ang{\al,k})=\ang{\al,k+1}$ for the given 
$\al$ and each $k\in\om,$ and $h(\ang{\beta,k})=\ang{\beta,k}$ 
whenever $\beta\not=\al.$ (Then $h(\i)=\i'$.) Thus $h$ shifts only 
the \dd\al th copy of $\dZ$ in $\bI$ but does not move anything else. 

The $h$ generates an order automorphism 
$Z\,\longmapsto\,Z':\peim\,\hbox{ onto }\,\peim$ in obvious way. 
We observe that $Y'=Y$ because $\supp Y=\eta\sq\xi_\al=\al\ti\dZ$.
 
We set $G'=\ans{Z':Z\in G}.$ Then $Y\in G',$ $G'$ is \dd\peim 
generic over $\gM,$ and moreover, $\gM[G']=\gM[G]$ because $h\in\gM.$ 

Let $\x'=\ang{\a'_\j:\j\in\bI}$ be defined from $G'$ as $\x$ was 
defined from $G.$ Then $\a'_\j=\a_{h(\j)}$ for all $\j;$ in 
particular \ (a) $\x'\res\eta=\x\res\eta\,,$ \ and \ (b) 
$\x'\res\xi'$ is a shift of $\x\res\xi,$ so that 
$\rL[\x'\res\xi']=\rL[\x\res\xi]$ in $\gN=\gM[G]=\gM[G']$. 

It follows from (a) and the choice of $Y$ that 
$\neg\;\psi(\x\res\eta,x)$ holds in $\gN$ provided a real $x$ 
satisfies $\rL[x]=\rL[\x'\res\xi']$ in $\gN.$ On the other hand, we 
have already got $\psi(\x\res\eta,x)$ in $\gN$ provided 
$\rL[x]=\rL[\x\res\xi]$ holds in $\gN.$ These two statements 
contradict each other by (b). 

This ends the proof of the ``neither'' part of 
Theorem~\ref{tge}.\qed\vspace{4mm}

This also ends the proof of Theorem~\ref{ge}.\qed


\begin{thebibliography}{99}
\label{pgl}

\bibitem{bl} J.~E.~Baumgartner and R.~Laver. Iterated perfect set 
forcing. {\it Annals of Mathematical Logic\/} 1979, 17, pp. 
271 -- 288.

\bibitem{g94} M.~Groszek. $\om_1^\ast$ as an initial segment 
of the c-degrees. {\it Journal of Symbolic Logic\/} 1994, 59, 
no 3, pp. 956 --976.

\bibitem{g} M.~Groszek. Applications of iterated perfect set 
forcing. {\it Annals of Pure and Applied Logic\/} 1988, 39, pp. 
19 -- 53.

\bibitem{gj} M.~Groszek and T.~Jech. Generalized iteration of 
forcing. {\it Transactions of the American Mathematical Society\/} 
1991, vol. 324, pp. 1 -- 26.

\bibitem{hkl} L.~A.~Harrington, A.~S.~Kechris, A.~Louveau. 
A Glimm -- Effros dichotomy for Borel equivalence relations. 
{\it J. Amer. Math. Soc.\/} 1990, 3, no 4, p. 903 --928.

\bibitem{h-det} G.~Hjorth. {\em A dichotomy for the definable 
universe\/}. (Preprint.)

\bibitem{hk} G.~Hjorth and A.~S.~Kechris. {\em Analytic 
equivalence relations and Ulm--type classification}.  
Department of Mathematics, Caltech.

\bibitem{k-sm} V.~Kanovei. {\it On a Glimm -- Effros dichotomy and 
an Ulm--type classification in Solovay model\/}. Logic Eprints, 
July 1995.

\bibitem{k-s11} V.~Kanovei. {\it On a Glimm -- Effros dichotomy 
theorem for Souslin relations in generic universes\/}. 
Logic Eprints, August 1995.

\bibitem{sa} G.~E.~Sacks. Forcing with perfect closed sets. 
{\it Axiomatic set theory, Part 1\/}. (Proceedings of Symposia in 
Pure mathematics, vol. 13. AMS, Providence, RI.) 1971, pp. 331 -- 
355. 

\bibitem{st} R.~M.~Solovay and S.~Tennenbaum. Iterated Cohen 
extensions and Souslin's problem. {\it Annals of Mathematics\/} 
1971, 94, pp. 201 -- 245.

\end{thebibliography}
\end{document}